\documentstyle[txmac,a4,%leqno,%
amssymb,%
case,%
twoside,%
pifont,%
nocaphead,%
epsf,%
mypic,times,mathptm]{article}

\advance\oddsidemargin by -1.9cm
\advance\evensidemargin by -1.9cm
\advance\textwidth by 3.8cm
\def\mynewtheo#1#2{%
\newtheorem{@#1}{#2}[section]%
\newenvironment{#1}{\begin{@#1}\rm}{\end{@#1}}}

\mynewtheo{lemma}{Lemma}
\mynewtheo{exer}{Exercise}
\mynewtheo{theo}{Theorem}
\mynewtheo{rem}{Remark}
\mynewtheo{defi}{Definition}
\mynewtheo{conj}{Conjecture}
\mynewtheo{corr}{Corollary}
\mynewtheo{prop}{Proposition}
\mynewtheo{question}{Question}
\mynewtheo{exam}{Example}

\newenvironment{theorem}{\begin{theo}}{\end{theo}}

\newenvironment{propo}{\begin{prop}}{\end{prop}}

\newenvironment{eqn}{\begin{equation}}{\end{equation}}

\begin{document}

\makeatletter

% eqlabel=(#1)
\newenvironment{myeqn*}[1]{\begingroup\def\@eqnnum{\reset@font\rm#1}%
\xdef\@tempk{\arabic{equation}}\begin{equation}\edef\@currentlabel{#1}}
{\end{equation}\endgroup\setcounter{equation}{\@tempk}\ignorespaces}

% eqlabel=#1
\newenvironment{myeqn}[1]{\begingroup\let\eq@num\@eqnnum
\def\@eqnnum{\bgroup\let\r@fn\normalcolor % an extremely UGLY hack !!!
\def\normalcolor####1(####2){\r@fn####1#1}%
%\show\reset@font
\eq@num\egroup}%
\xdef\@tempk{\arabic{equation}}\begin{equation}\edef\@currentlabel{#1}}
{\end{equation}\endgroup\setcounter{equation}{\@tempk}\ignorespaces}

% eqlabel=(eqnnr) \qed
\newenvironment{myeqn**}{\begin{myeqn}{(\arabic{equation})\es\es\mbox{\qed}}\edef\@currentlabel{\arabic{equation}}}
{\end{myeqn}\stepcounter{equation}}

\newcommand{\mybin}[2]{\text{$\Bigl(\begin{array}{@{}c@{}}#1\\#2%
\end{array}\Bigr)$}}
\newcommand{\mybinn}[2]{\text{$\biggl(\begin{array}{@{}c@{}}%
#1\\#2\end{array}\biggr)$}}

\def\overtwo#1{\mbox{\small$\mybin{#1}{2}$}}
\newcommand{\mybr}[2]{\text{$\Bigl\lfloor\mbox{%
\small$\displaystyle\frac{#1}{#2}$}\Bigr\rfloor$}}
\def\mybrtwo#1{\mbox{\mybr{#1}{2}}}

\def\myfrac#1#2{\raisebox{0.2em}{\small$#1$}\!/\!\raisebox{-0.2em}{\small$#2$}}

\author{A. Stoimenow\footnotemark[1]\\[2mm]
\small Department of Mathematics, \\
\small University of Toronto,\\
\small Canada M5S 3G3\\
\small e-mail: {\tt stoimeno@math.toronto.edu}\\
\small WWW: {\hbox{\tt http://www.math.toronto.edu/stoimeno/}}
}

{\def\thefootnote{\fnsymbol{footnote}}
\footnotetext[1]{Supported by a DFG postdoc grant.}
}

\title{%
\large\bf \uppercase{Square numbers, spanning trees and}\\[2mm]
\uppercase{invariants of achiral knots}\\[3mm]
{\small\it This is a preprint. I would be grateful
for any comments and corrections!}}

\date{\large Current version: \curv\ \ \ First version:
\makedate{10}{1}{2000}}

\maketitle

\makeatletter

\let\vn\varnothing
\let\point\pt
\let\ay\asymp
\let\pa\partial
\let\al\alpha
\let\be\beta
\let\Dl\Delta
\let\Gm\Gamma
\let\gm\gamma
\let\de\delta
\let\dl\delta
\let\eps\epsilon
\let\lm\lambda
\let\Lm\Lambda
\let\sg\sigma
\let\vp\varphi
\let\om\omega
\let\diagram\diag

\let\sm\setminus
\let\tl\tilde
\let\wt\widetilde
\def\ncap{\not\mathrel{\cap}}
\def\sgn{\text{\rm sgn}\,}
\def\spn{\mathop {\operator@font span}}
\def\md{\max\deg}
\def\mc{\max\cf}
\def\mcf{\min\cf}
\def\Lra{\Longrightarrow}
\def\lra{\longrightarrow}
\def\so{\Rightarrow}
\def\So{\Longrightarrow}
\def\nin{\not\in}
\let\ds\displaystyle
\let\llra\longleftrightarrow
\let\reference\ref
\let\dt\det

\long\def\@makecaption#1#2{%
   % \tm
   \vskip \abovecaptionskip 
   {\let\label\@gobble
   \let\ignorespaces\@empty
   \xdef\@tempt{#2}%
   }%
   \ea\@ifempty\ea{\@tempt}{%
   \sbox\@tempboxa{%
      \fignr#1#2}%
      }{%
   \sbox\@tempboxa{%
      {\fignr#1:}\capt\ #2}%
      }%
   \ifdim \wd\@tempboxa >\captionwidth {%
      \centerline{\parbox{\captionwidth}{\box \@tempboxa}}%
      %\rightskip=\@captionmargin\leftskip=\@captionmargin
      %\unhbox\@tempboxa\par
     }%
   \else
      \centerline{\box \@tempboxa}%
      % \hbox to\captionwidth{\hfil\box\@tempboxa\hfil}%
   \fi
   \vskip \belowcaptionskip
   }%
\def\fignr{\small\sffamily\bfseries}%
\def\capt{\small\sffamily}%

% \long\def\@makecaption#1#2{%
%    % \tm
%    \vskip 10pt
%    {\let\label\@gobble
%    \let\ignorespaces\@empty
%    \xdef\@tempt{#2}%
%    %\typeout{`#2'}%
%    }%
%    \ea\@ifempty\ea{\@tempt}{%
%    \setbox\@tempboxa\hbox{%
%       \fignr#1#2}%
%       }{%
%    \setbox\@tempboxa\hbox{%
%       {\fignr#1:}\capt\ #2}%
%       }%
%    \ifdim \wd\@tempboxa >\captionwidth {%
%       \rightskip=\@captionmargin\leftskip=\@captionmargin
%       \unhbox\@tempboxa\par}%
%    \else
%       \hbox to\captionwidth{\hfil\box\@tempboxa\hfil}%
%    \fi}%
% %
% \def\fignr{\small\sffamily\bfseries}%
% \def\capt{\small\sffamily}%

\newdimen\@captionmargin\@captionmargin2cm\relax
\newdimen\captionwidth\captionwidth0.8\hsize\relax

\def\eqref#1{(\protect\ref{#1})}

\def\proof{\@ifnextchar[{\@proof}{\@proof[\unskip]}}
\def\@proof[#1]{\noindent{\bf Proof #1.}\enspace}

\def\hint{\noindent Hint: }
\def\problem{\noindent{\bf Problem.} }

\def\@mt#1{\ifmmode#1\else$#1$\fi}
\def\qed{\hfill\@mt{\Box}}
\def\qqed{\hfill\@mt{\Box\enspace\Box}}

\def\cI{{\cal I}}
\def\cU{{\cal U}}
\def\cC{{\cal C}}
\def\cP{{\cal P}}
\def\cL{{\cal L}}
\def\fg{{\frak g}}
\def\Kr{\mathop{\operator@font Kr}}
\def\cf{\mathop{\operator@font cf}}
\def\cZ{{\cal Z}}
\def\cD{{\cal D}}
\def\bQ{{\Bbb Q}}
\def\bR{{\Bbb R}}
\def\cE{{\cal E}}
\def\bZ{{\Bbb Z}}
\def\bN{{\Bbb N}}

\def\A{\diag{1em}{2}{1.5}{\picline{0 0}{2 0}\picline{0 1.5}{2 1.5}}}
\def\B{\diag{1em}{2}{1.5}{\picellipsearc{0 0.75}{0.8 0.75}{-90 90}
\picellipsearc{2 0.75}{0.8 0.75}{90 -90}}}
\def\AT{\diag{0.7em}{1.5}{2}{\picline{0 0}{0 2}\picline{1.5 0}{1.5 2}}}
\def\BT{\diag{0.7em}{1.5}{2}{\picellipsearc{0.75 0}{0.75 0.8}{0 180}
\picellipsearc{0.75 2}{0.75 0.8}{180 0}}}

\def\geni{\diag{3mm}{2}{2}{
  \picline{0 0}{0 2}
  \picline{1 0}{1 2}
  \picline{2 0}{2 2}
}}
\def\genii{
\diag{3mm}{2}{2}{
  \picline{0 0}{0 2}
  \piccirclearc{1.5 0}{0.5}{0 180}
  \piccirclearc{1.5 2}{0.5}{0 180 x}
}}
\def\geniii{
\diag{3mm}{2}{2}{
  \piccirclearc{0.5 0}{0.5}{0 180}
  \piccirclearc{0.5 2}{0.5}{0 180 x}
  \picline{2 0}{2 2}
}}
\def\geniv{
\diag{3mm}{2}{2}{
  \piccirclearc{1.5 0}{0.5}{0 180}
  \piccirclearc{0.5 2}{0.5}{0 180 x}
  \picline{0 0}{2 2}
}}
\def\genv{
\diag{3mm}{2}{2}{
  \piccirclearc{0.5 0}{0.5}{0 180}
  \piccirclearc{1.5 2}{0.5}{0 180 x}
  \picline{0 2}{2 0}
}}
\def\rb{\ry{7mm}\raisebox{0.5em}}
\def\xd#1{\mbox{\ding{#1}}}

\def\bysame{\same[\kern2cm]\,}

\def\br#1{\left\lfloor#1\right\rfloor}
\def\BR#1{\left\lceil#1\right\rceil}

\def\abstractname{}

\@addtoreset {footnote}{page}

\renewcommand{\section}{%
   \@startsection
         {section}{1}{\z@}{-1.5ex \@plus -1ex \@minus -.2ex}%
               {1ex \@plus.2ex}{\large\bf}%
}
\renewcommand{\@seccntformat}[1]{\csname the#1\endcsname .
\quad}

\def\bC{{\Bbb C}}
\def\bP{{\Bbb P}}

\def\epsfs#1#2{{\catcode`\_=11\relax\ifautoepsf\unitxsize#1\relax\else
\epsfxsize#1\relax\fi\epsffile{#2.eps}}}
\def\epsfsv#1#2{{\vcbox{\epsfs{#1}{#2}}}}
\def\vcbox#1{\setbox\@tempboxa=\hbox{#1}\parbox{\wd\@tempboxa}{\box
  \@tempboxa}}
\def\p{\epsfsv{2cm}}

\def\@test#1#2#3#4{%
  \let\@tempa\go@
  \@tempdima#1\relax\@tempdimb#3\@tempdima\relax\@tempdima#4\unitxsize\relax
  \ifdim \@tempdimb>\z@\relax
    \ifdim \@tempdimb<#2%
      \def\@tempa{\@test{#1}{#2}}%
    \fi
  \fi
  \@tempa
}

\def\go@#1\@end{}
\newdimen\unitxsize
\newif\ifautoepsf\autoepsffalse

\unitxsize4cm\relax
\def\epsfsize#1#2{\epsfxsize\relax\ifautoepsf
  {\@test{#1}{#2}{0.1 }{4   }
		{0.2 }{3   }
		{0.3 }{2   }
		{0.4 }{1.7 }
		{0.5 }{1.5 }
		{0.6 }{1.4 }
		{0.7 }{1.3 }
		{0.8 }{1.2 }
		{0.9 }{1.1 }
		{1.1 }{1.  }
		{1.2 }{0.9 }
		{1.4 }{0.8 }
		{1.6 }{0.75}
		{2.  }{0.7 }
		{2.25}{0.6 }
		{3   }{0.55}
		{5   }{0.5 }
		{10  }{0.33}
		{-1  }{0.25}\@end
		\ea}\ea\epsfxsize\the\@tempdima\relax
		\fi
		}

\let\old@tl\~\def\~{\raisebox{-0.8ex}{\tt\old@tl{}}}
\let\lra\longrightarrow
\let\sm\setminus
\let\eps\varepsilon
\let\ex\exists
\let\fa\forall
\let\ps\supset

\def\rs#1{\raisebox{-0.4em}{$\big|_{#1}$}}

{\let\@noitemerr\relax
\vskip-2.7em\kern0pt\begin{abstract}
\noindent{\bf Abstract.}\enspace
We give constructions to realize an odd number, which is representable
as sum of two squares, as determinant of an achiral knot, thus
proving that these are exactly the numbers occurring as such
determinants. Later we study which numbers occur as determinants
of prime alternating achiral knots, and obtain a complete
result for perfect squares.
% 
% We examine and partially confirm some questions on  properties of the
% the Alexander and HOMFLY polynomial of achiral knots. In particular we
% show that determinants of achiral knots are exactly the odd numbers
% representable as sums of two squares.
Using the checkerboard coloring, then an application is given to the
number of spanning trees in planar self-dual graphs. Another application
are some enumeration results on achiral rational knots. Finally, we
describe the leading coefficients of the Alexander and skein
polynomial of alternating achiral knots.\\[1mm]
\end{abstract}
}\vspace{7mm}

{\parskip0.2mm\tableofcontents}
\vspace{7mm}

\section{Introduction}

The main problem of knot theory is to distinguish knots (or links),
i.e., smooth embeddings of $S^1$ (or several copies of it)
into $\bR^3$ or $S^3$ up to isotopy. A main tool
for this is to find \em{invariants} of knots, i.e., maps of
knot diagrams into some algebraic structure, which are invariant
under Reidemeister's moves. A family of most popular such invariants
are the polynomial invariants, associating to each knot an element
in some one- or two-variable (Laurent) polynomial ring over $\bZ$.
Given a knot invariant, beside distinguishing knots with it, one
is also interested which properties of knots it reflects, and in
which way.

One of the most intuitive ways to associate to a knot (or link) another
one is to consider its \em{obverse}, or mirror image, obtained by
reversing the orientation of the ambient space. The knot (or link)
is called \em{achiral} (or synonymously \em{amphicheiral}),
if it coincides (up to isotopy) with its mirror image, and \em{chiral}
otherwise. When considering orientation of the \em{knot}, then we
distinguish among achiral knots between $+$achiral and $-$achiral ones,
dependingly on whether the deformation into the mirror image preserves
or reverses the orientation of the knot. (For links one has to
attach a sign to each component, i.e. embedded circle, and take
into account possible permutations of the components.)

When the \em{Jones polynomial} $V$ \cite{Jones} appeared in 1984, one of
its (at that time) spectacular features was that it was (in general)
able to distinguish between a knot and its obverse by virtue of having
distinct values on both, and (hence) so were its generalizations,
the \em{HOMFLY, or skein, polynomial} $P$ \cite{HOMFLY}
% (we use henceforth the convention of \cite{LickMil} for it)
and the \em{Kauffman polynomial} $F$ \cite{Kauffman}.
The $V$, $P$ and $F$ polynomials of achiral knots have
the special property to be \em{self-conjugate}, that is,
invariant when one of the variables is replaced by its inverse.
Their decades-old predecessor, the \em{Alexander polynomial} $\Dl$
\cite{Alexander}, a knot invariant with values in $\bZ[t,t^{-1}]$,
was known always to take the same value on a knot and its mirror image.
Nevertheless, contrarily to the common belief, $\Dl$ can also be
used to detect chirality (the property of a knot to be distinct from its
mirror image) by considering its value $\Dl(-1)$,
called \em{determinant}. 

The aim of this paper is to study invariants of achiral knots and
to relate some properties of their determinants to the classical topic
in number theory of representations of integers as sums of two squares.

In \S\reference{S2} we begin with recalling a criterion for the
Alexander polynomial of an achiral knot via the determinant,
which follows from Murasugi's work on the signature and the
Lickorish-Millett value of the Jones polynomial. This conditions
show that, paradoxly formulated, although the Alexander polynomial
cannot distinguish between a knot and its mirror image, it can still
sometimes show that they are distinct.

After collecting
some number theoretic preliminaries in \S\reference{Nth}, we show
then in \S\reference{S31} that the condition of \S\reference{S2}
is in fact a reduction modulo 36 of the exact arithmetic description
of numbers, occurring as determinants of achiral knots. Namely, an
odd natural number is the determinant of an achiral knot if and only
if it is the sum of two squares. The `only if' part of this
statement was an observation of Hartley and Kawauchi in \cite{HarKaw}.
Our aim will be to show the `if' part, that is, given a sum of two
squares, to realize it as the determinant of an achiral knot (theorem
\reference{Th1}). The main tool used is the definition of
the determinant by means of Kauffman's state model for
the Jones polynomial \cite{Kauffman2}.

Then we attempt to refine our construction, by
producing achiral knots (of given determinant) with additional
properties: prime and/or alternating. Although it turns out,
that one of these properties can always easily be achieved, the
situation reveals much harder when demanding them both altogether.
We investigate this problem in \S\reference{Spa}. Now, the
correspondence of \S\reference{S31} does not hold completely,
and there are exceptional values of the determinant, that cannot
be realized. To show that 9 and 49 are such, we prove a quadratic
improvement of Crowell's (lower) bound for the determinant of an
alternating knot in terms of its crossing number \cite{Crowell2},
in the case the knot is achiral (proposition \reference{prp2}).
We obtain then a complete result about which perfect squares can
be realized as determinants of prime alternating achiral knots
(theorem \reference{Thsq}).

In \S\reference{Su1} we consider the problem to describe determinants
of unknotting number one achiral knots. In this case the description
is even less clear, as we show by several examples.

Then we give some applications, including enumeration
results on rational knots in \S\reference{Senum},
and a translation of the previously established properties
to the number of spanning trees in planar self-dual graphs
in \S\reference{par5}.

% In \S\reference{S3}
% 
% and develop in further in
% \S\reference{S3} to an exact arithmetic description of determinants
% of achiral knots, providing the reverse direction to an observation
% of Hartley and Kawauchi in \cite{HarKaw}. In \S\reference {par5} the
% application to spanning trees in planar graphs is discussed,
% and some partial results
% on generalizations are given in \S\reference{par4}.

Subsequently, in \S\reference{S4} we prove some further (at parts still
remaining conjectural) properties of the Alexander and skein polynomial
of at least large classes of achiral knots, which would allow to decide
about chirality (the lack of an isotopy to the mirror image)
in a yet different way, at least for these knot classes. These
properties concern the leading coefficients of the polynomials,
and are closely related to Murasugi's $*$-product. We prove
in particular that perfect squares are exactly the numbers occurring
as leading coefficients of the Alexander polynomial of alternating
achiral knots, thus improving the previously known necessary
condition of non-primeness (corollary \reference{crx}).
These results have been obtained with the same arguments
independently (but somewhat later) by C. Weber and
Q.\ H.\ C\^am V\^an \cite{VW}.

Several open problems are suggested during the discussion throughout
the paper. These problems appear to be involved enough already for
knots, so that we waived on an analogous study of links (which
are the cases covering the even natural numbers). For links,
also the unpleasant issue of component orientations becomes relevant.

% previous work, mainly of Cromwell
% \cite{Cromwell} and Murasugi-Przytycki \cite{MurPrz}, but apparently
% have not been drawn attention to explicitly before.
% We will make some remarks how the conjectured properties follow for
% rational (and some other) knots. However,
% % but will not put central emphasis of these results,
% the most significant open part of the problems appears to deserve
% more attention than the solution for these partial cases.

\section{Detecting chirality with the Alexander polynomial\label{S2}}

In the following we will be concerned with the value $\Dl(-1)$
of the Alexander polynomial, where $\Dl$ is normalized
so that $\Dl(t)=\Dl(1/t)$ and $\Dl(1)=1$. Up to sign, this numerical
invariant can be interpreted as the order of the homology
group (over $\bZ$) of the (double) branched covering of $S^3$ over $K$
associated to the canonical homomorphism $\pi_1(S^3\sm K)\to\bZ_2$ and
carries the name ``\em{determinant}'' because of its expression (up to
sign) as the determinant of a Seifert \cite[p.\ 213]{Rolfsen} or Goeritz
\cite{Goeritz} matrix.

To introduce some mathematical notations of the objects thus
occurring, let $D_K$ be the double branched cover of $S^3$ over a knot
$K$, and let $H_1(D_K)$ be its homology group (over $\bZ$). 
We write then $\dt(K)=|\,\Dl_K(-1)\,|=|\,H_1(D_K)\,|$.

We start first by description of two special cases of the exact
property of the determinant of achiral knots, which we will formulate
subsequently, because they have occurred in independent contexts and
deserve mention in their own right. They allow to decide
about chirality of a knot $K$, at least for
$\myfrac{11}{18}$ of the possible values of $\Dl_K(-1)$.

There is an observation (originally likely, at least implicitly,
due to Murasugi \cite{Murasugi}, and applied explicitly
in \cite{nvi}), using the sign of the value $\Dl(-1)$
(with $\Dl$ normalized as said). The information of this sign is
equivalent to the residue $\sg\bmod 4$, where $\sg$ denotes
the \em{signature}. Whenever $\Dl(-1)<0$, we have $\sg\equiv 2\bmod 4$,
so in particular $\sg\ne 0$, and the knot cannot be achiral.
This argument works e.g. for the knot $9_{42}$ in the tables
of \cite[appendix]{Rolfsen}, which became famous
by sharing the same $V$, $P$ and $F$ polynomial with its obverse,
since its polynomials are all self-conjugate.

Another way to deduce chirality from the determinant 
is to use the sign of the Lickorish-Millett value
$V\bigl(e^{\pi i/3}\bigr)$ \cite{LickMil2}. Attention to it
was drawn in \cite{Traczyk}, where it was used to calculate
unknotting numbers. Using some of the ideas there,
in \cite{unkn1} we observed that this sign
implies that if for an achiral knot $3\mid\Dl(-1)$, then
already $9\mid\Dl(-1)$. Thus for example also the chirality
of $7_7$ can be seen already from its Alexander polynomial,
as in this case $\Dl(-1)=21$ (although the Murasugi trick does not
work here, and indeed $\sg=0$).

Combining both criteria, we arrive in summary to

\begin{prop}\label{pR1}
For any achiral knot $K$ we have $\big|\Dl_K(-1)\big|
\bmod 36\in\{1,\ 5,\ 9,\ 13,\ 17,\ 25,\ 29\}$. \qed
\end{prop}

An easy verification shows that all these residues indeed occur.

In view of these opportunities to extract chirality information out of
$\Dl$, it appears appropriate to introduce a clear distinction between
the terms `detecting chirality of $K$', meant in the sense `showing
that $K$ and $!K$ are not the same knot' (which can be achieved
by the above tricks) and `distinguishing between $K$ and $!K$', meant
in the sense `identifying for a given diagram, known \em{a priori}
to belong to either $K$ or $!K$, to which one of both it belongs'
(what they cannot accomplish, but what is the usually imagined
situation where some of the other polynomials is not self-conjugate).

% achir.tex

Here is a small arithmetic consequence. It is elementary, but is
included because of its knot theoretical interpretation and as it is
the starting point of exhibiting some more interesting phenomena
described in the next sections.

Recall, that a knot $K$ is 
\em{rational} (or $2$-bridge), if it has an embedding with a
Morse function having only four critical points (2 maxima and 2 minima).
Such knots were classified by Schubert \cite{Schubert}, and can
be alternatively described by their Conway notation \cite{Conway}.
See for example \cite[\S 2.3]{Adams} for a detailed description.
It is well-known that rational knots are alternating (see
\cite[proposition 12.14, p.\ 189]{BurZie}).

\begin{corr}\label{corr2.1}
Let $p/q$ for $(p,q)=1$, $p$ odd be expressible as the continued
fraction
\[
[[a_1,\dots,a_n,a_n,\dots,a_1]]\,:=\,
a_1+\frac{1}{a_2+\frac{1}{\dots\, a_2+\frac{1}{a_1}}}
\]
for a palindromic sequence $(a_1,\dots,a_n,a_n,\dots,a_1)$ of even
length (with the usual conventions $\frac{1}{0}=\infty$, 
$\infty+*=\infty$, $\frac{1}{\infty}=0$ for the degenerate cases).
Then $|p|\equiv 1$ or $5\bmod 12$.
\end{corr}

\proof Observe that $|p|$ is the determinant of the achiral rational
knot with Conway notation $(a_1\dots a_na_n\dots a_1)$.
The above proposition \reference{pR1} leaves us only with
explaining why $9\nmid p$. The implication $3\mid\det(K)\Lra
9\mid\det(K)$ for $K$ achiral using $V\bigl(e^{\pi i/3}\bigr)$
follows from the fact that the number of torsion coefficients
divisible by $3$ of the $\bZ$-module $H_1(D_K)$, counted by
$V\bigl(e^{\pi i/3}\bigr)$, is even. However, for a rational
knot $K$, $H_1(D_K)$ is cyclic and non-trivial ($D_K$ is a lens
space), so that there is only one torsion number at all. Thus
$H_1(D_K)$ for any achiral rational knot $K$ cannot have any
$3$-torsion. \qed

\section{Number theoretic preliminaries\label{Nth}}

According to a claim of Fermat, written about 1640 on the margins of
his copy of Euclid's ``Elements'', proved in 1754 by Euler, and
further simplified to the length of ``one sentence'' in \cite{Zagier3},
% In 1754, Euler showed one of Fermat's claims written on the margins of
% his copy of Euclid's ``Elements'', namely, that
any prime of the form $4x+1$ can be written as the sum of two squares.
% Today, this result is known to be a
% special case of the more general statement that any natural number $n$
More generally, any natural number $n$
is the sum of two squares if and only if any prime of the form $4x+3$
occurs in the prime decomposition of $n$ with an even power, and
it is the sum of the squares of two coprime numbers if and only if
such primes do not occur at all in the prime decomposition of $n$.

The number of representations as the sum of two squares is given by
the formula
\begin{eqn}\label{*}
r_2(n)\,:=\,\frac{1}{4}\,\#\{\,(m_1,m_2)\in\bZ^2\,:\,m_1^2+m_2^2=n\,\}\,=\,
\sum_{d\mid n}\left( \frac{-4}{d}\right)\,=\, \#\{\,x\in\bN\,:
\,4x+1\mid n\,\}-\#\{\,x\in\bN\,:\,4x+3\mid n\,\}\,,
\end{eqn}
which has also an interpretation in the theory of modular forms
(see \cite[(16.9.2) and theorem 278, p.\ 275]{HardyWright} and
\cite{Zagier2}). 

A number theoretic explanation of \eqref{*} is as follows:
If we denote by $\chi$ the (primitive) character $\left(\frac
{-4}{\,\cdot\,}\right)$, with $-4$ being the discriminant of the
field of the Gau\ss{} numbers $\bQ[i]$, we have for $\Re(s)>1$,
using that $\bQ[i]$ has class number $1$ and 4 units, that
\begin{eqn}
\sum_{n=1}^{\infty}\frac{r_2(n)}{n^s} \, = \,
\zeta_{\bQ[i]}(s) \, = \, \zeta(s)\,L(s,\chi) 
\, = \, 
\prod_{\scbox{$p\ \mbox{prime}$}}\,\frac{1}{(1-p^{-s})(1-\chi(p)
p^{-s})}\,,
\end{eqn}
from which the formula follows by considering the Taylor expansion
in $p^{-s}$ of the different factors in the product.
(This series converges for $\Re(s)>1$.)

The $\zeta$-function identities also give a formula for
\begin{eqn}%\label{r20}
r_2^0(n)\,:=\,\frac{1}{4}\,
\#\{\,(a,b)\in\bZ^2\,:\,(a,b)=1,\, a^2+b^2=n\,\},
\end{eqn} 
the number of representations of $n$ as the sum of squares of
coprime numbers.

We have
\begin{eqn}
\sum_{n=1}^{\infty}\frac{r_2^0(n)}{n^s} \, = \,
\frac{\zeta_{\bQ[i]}(s)}{\zeta(2s)} 
\, = \, \frac{\zeta(s)\,L(s,\chi)}{\zeta(2s)} 
\, = \, 
\prod_{\scbox{$p\ \mbox{prime}$}}\,\frac{1-p^{-2s}}{(1-p^{-s})(1-\chi(p)
p^{-s})}\, = \, (1+2^{-s})
\prod_{\scbox{$\begin{array}{c}p\equiv 1\,(4)\\\mbox{prime}
\end{array}$}} \,\frac{1+p^{-s}}{1-p^{-s}}\,.
\end{eqn}
Thus
\begin{eqn}\label{r02}
r_2^0(n) \, = \,% \nonumber \qquad =\,
\left\{\begin{array}{ll}
2^k & \mbox{if } n=p_1^{k_1}\cdot\dots\cdot p_k^{d_k} \mbox{ or }
2p_1^{k_1}\cdot\dots\cdot p_k^{d_k} \mbox{ with $p_1<p_2<\dots<p_k$
primes $\equiv 1\bmod 4$ and $d_i>0$} \\
0 & \mbox{else}
\end{array}\right..%\}
\end{eqn}
Note, that for $n>1$,
% Similarly
% one deduces \eqref{*} from
% \[
% \sum_{n=1}^{\infty}\frac{r_2(n)}{n^s} \, = \,
% \zeta_{\bQ[i]}(s)\, = \, \zeta(s)\,L(s,\chi)\,.
% \]
% 
% Define
\begin{eqn}\label{r20}
r_2^0(n)\,=\,
\#\{\,(a,b)\in\bN^2\,:\,(a,b)=1,\, a^2+b^2=n\,\}\,,
\end{eqn} 
% to be the number of representations of $n$ as the sum of squares of
% coprime numbers.
and 
% A formula for $r_2^0(n)$ is known from elementary
% number theory. See \cite[theorem 367, \S 20.3]{HardyWright} for a
% partial statement. The formula follows for prime powers by a simple 
% inclusion-exclusion argument from \eqref{*}, and for arbitrary numbers
% by multiplicativity. % (the number theoretic formulation of this
% % argument given below is due to D.\ Zagier).
% 
% \begin{theo}
for $n>2$ we have
\begin{eqn}
\frac{r_2^0(n)}2\,=\,
\#\{\,(a,b)\in\bN^2\,:\,(a,b)=1,\, a\le b,\, a^2+b^2=n\,\}\,.
\end{eqn}
% \end{theo}

% These numbers can be given~-- to the number theorist's eyes most
% pleasantly~-- by the Euler product form of the associated $L$-series:
% \[
% 1+2\sum_{n=1}^\infty\,\frac{c_n}{n^s}\quad =\,%\frac{1}{2}\Bigl[
% \,\prod_{\scbox{$\begin{array}{c}p\equiv 1\,(4)\\\mbox{prime}
% \end{array}$}} \,\frac{1+p^{-s}}{1-p^{-s}}\,.
% % \Bigr]\,.
% \]
% (This series converges for $\Re(s)>1$.)
% A number theoretic interpretation of this fact is as follows:
% If we denote by $\chi$ the (primitive) character $\left(\frac
% {-4}{\,\cdot\,}\right)$, with $-4$ being the discriminant of the
% field of the Gau\ss{} numbers $\bQ[i]$, we have for $\Re(s)>1$,
% using that $\bQ[i]$ has class number $1$ and 4 units, that
% \begin{eqn}
% \sum_{n=1}^{\infty}\frac{r_2^0(n)}{n^s} \, = \,
% \frac{\zeta_{\bQ[i]}(s)}{\zeta(2s)} 
% \, = \, \frac{\zeta(s)\,L(s,\chi)}{\zeta(2s)} 
% \, = \, 
% \prod_{\scbox{$p\ \mbox{prime}$}}\,\frac{1-p^{-2s}}{(1-p^{-s})(1-\chi(p)
% p^{-s})}\,,
% \end{eqn}
% from which the formula follows by considering the Taylor expansion
% in $p^{-s}$ of the different factors in the product. Similarly
% one deduces \eqref{*} from
% \[
% \sum_{n=1}^{\infty}\frac{r_2(n)}{n^s} \, = \,
% \zeta_{\bQ[i]}(s)\, = \, \zeta(s)\,L(s,\chi)\,.
% \]

\section{Sums of two squares and determinants of achiral knots\label{S3}}

\subsection{Realizing sums of two squares as determinants\label{S31}}

The aim of this section is to establish a, partially conjectural,
correspondence between sums of two squares and the determinant
of achiral links. The study of this relation 
was first initiated in \cite{HarKaw}, where it was observed
that a result of Goeritz \cite{Goeritz} implies that the determinant
of an achiral knot is the sum of two squares. We shall here
show the converse. In fact, we have

\begin{theo}\label{Th1}
An odd natural number $n$ occurs as determinant of an achiral knot $K$
if and only if $n$ is the sum of two squares $a^2+b^2$. More
specifically,
% if $n$ is the sum of two squares $a^2+b^2$, then $n$ is the determinant
% of an achiral knot $K$. In fact,
% if $a$ and $b$ can be chosen to be both non-zero, then
$K$ can be chosen to be alternating or prime, and if one can choose
$a$ and $b$ to be coprime, then we can even take
$K$ to be rational (or $2$-bridge).
\end{theo}

Note, that we have given another
argument from that of Hartley and Kawauchi for the reverse implication
``modulo 36'': it was observed above how the signature
and the Lickorish-Millett value of the Jones polynomial imply that
if $n$ is the determinant of an achiral knot, then $n\bmod 36\in\{1,\ 
5,\ 9,\ 13,\ 17,\ 25,\ 29\}$. These are exactly the congruences 
which odd sums of two squares leave modulo 36. Clearly, not every number
satisfying these congruences is the sum of two squares. The simplest
example is $77$. And indeed, this number does not occur as determinant
of any achiral knot of $\le 16$ crossings. % In fact, no such knot
% provided a counterexample to conjecture \ref{cj2}, so that there is
% also much empirical evidence for it.

% \begin{corr}\label{cOr1}
% Determinants of achiral knots are exactly the odd numbers
% representable as sums of two squares.
% \end{corr}

% We will in the following restrict our attention to knots and odd $n$
% and prove one special case of our conjecture \ref{cj1}.

Additionally to the general case, we also have a complete statement for
rational knots. (An analogue for arborescent knots seems possible by
applying the classification result of Bonahon and Siebenmann
\cite{BonSie}.)

\begin{theo}\label{Th2}
An odd natural number $n$ is the determinant of an achiral rational knot
if and only if it is the sum of the squares of two coprime numbers.
\end{theo}

The coprimality condition is clearly restrictive~-- for example $49$
and $121$ are not sums of the squares of two coprime numbers. Moreover,
it also implies the congruence modulo 12 proved in corollary
\reference{corr2.1}.

Fermat's theorem can be now knot-theoretically reformulated for example
as 

\begin{corr}
If $n=4x+1$ is a prime, then there is a rational achiral knot with 
determinant $n$.
\end{corr}

\proof We have $n=a^2+b^2$ and as $n$ is prime, $a$ and $b$
must be coprime. \qed

We start by a proof of theorem \ref{Th2}.
For this recall  Krebes's invariant defined in \cite{Krebes}.
Any tangle $T$ can be expressed by its coefficients in the
the Kauffman bracket skein module of the room with four
in/outputs (see \cite{gwg}):
\begin{eqn}\label{ABC}
\diag{5mm}{1}{1.5}{
  \piclinewidth{40}
  \pictranslate{0.5 0.75}{
    \picline{-0.4 -0.75}{-0.4 0.75}
    \picline{0.4 -0.75}{0.4 0.75}
    \picfilledbox{0 0}{1 1}{$T$}
  }
}
\quad=\quad A\,\,\AT\enspace +\enspace B\,\,\BT\,.
\end{eqn}

\begin{defi}
For a tangle $
\diag{4mm}{2.4}{2.4}{
  \pictranslate{1.2 1.2}{
  \picmultigraphics[rt]{4}{90}{
  \picline{1.0 45 polar}{1.6 45 polar}
  }
  \picfilledbox{0 0}{1.6 d}{$T$}
}}
$, we call 
$\diag{4mm}{2.4}{3}{
  \picmultigraphics{2}{0 2}{
    \picellipse{1.2 0.5}{1.2 0.5}{}
  }
  \picfilledbox{1.2 1.5}{1.6 d}{$T$}
}$ the \em{numerator closure} of $T$ and
$\diag{4mm}{3}{2.4}{
  \picmultigraphics{2}{0 2 x}{
    \picellipse{1.2 0.5 x}{1.2 0.5 x}{}
  }
  \picfilledbox{1.2 1.5 x}{1.6 d}{$T$}
}$ the \em{denominator closure} of $T$.
\end{defi}

When specializing the bracket variable to $\sqrt{i}$ ($i$ denotes
henceforth $\sqrt{-1}$), $A$ and $B$ in \eqref{ABC} become scalars.

Then Krebes's invariant can be defined by
\[
\Kr(T)\,:=\,\frac{A}{B}\,=\,(A,B)\in \bZ\times\bZ/(p,q)\sim(-p,-q)\,.
\]
The denominator and numerator of this ``fraction'' give the
determinants of the two closures of $T$.

\proof[of theorem \reference{Th2}]
A rational achiral knot $(a_1\dots a_n a_n \dots a_1)$
is of the form 
\begin{eqn}\label{Tsum}
K\quad=\quad
\diag{6mm}{6}{3}{
  \pictranslate{3 1.5}{
       \picmultigraphics[S]{2}{1 -1}{
           \picmultiline{0.12 1 -1.0 0}{2 -0.5}{1 0.5}
           \picmultigraphics[S]{2}{-1 1}{
                \picellipsearc{-2 -1.0}{1 0.5}{90 270}
           }
           \picline{-2 -1.5}{2 -1.5}
           \picline{-2 -0.5}{1 -0.5}
      }
  }
  \pictranslate{1.7 1.5}{
    \picrotate{-90}{
      \picfilledbox{0 0}{1.4 1.4}{$\overline T$}
    }
  }
  \picfilledbox{4.3 1.5}{1.4 1.4}{$T$}
}\quad,
\end{eqn}
where $T=(a_1 \dots a_n)$ is a rational tangle and $\overline T$ its
mirror image. Because of connectivity
reasons $T$ must be of homotopy type 
$\,\diag{3mm}{2}{2}{
  \piccirclearc{-1 1}{1.41}{-45 45}
  \piccirclearc{3 1}{1.41}{135 -135}
}$ or $
\diag{3mm}{2}{2}{
  \piccirclearc{1 -1}{1.41}{45 135}
  \piccirclearc{1 3}{1.41}{-135 -45}
}$, i.e. $\Kr(T)=
[[a_1,\dots,a_n]]=\frac{p}{q}$ with $(p,q)=1$ and exactly one of $p$ and
$q$ is odd. Thus $K$ is the numerator closure of
the tangle sum $T+\overline{T}^T$, where ${\,.\,}^T$ denotes
transposition. By the calculus introduced by Krebes, his
invariant is additive under tangle sum, and invertive under
transposition, and so we have 
\[
\frac{\det(K)}{*}\,=\,\Kr(T+\overline{T}^T)=\Kr(T)+\frac{1}{\Kr(T)}
=\frac{p}{q}+\frac{q}{p}=\frac{p^2+q^2}{pq}\,.
\]
To justify our choice of sign in this calculation, that is, that the
% (Admittedly, we did not take care of the signs, but to see
determinant is $p^2+q^2$ rather than $p^2-q^2$, it suffices to keep in
mind that the diagram \eqref{Tsum} is alternating and in calculating
the bracket of alternating diagrams no cancellations occur, as
explained also in \cite{Krebes}. Thus we have the `only if' part.

For the `if' part note that if $a$ and $b$ are coprime, then $\frac{a}
{b}$ can be expressed by an continued fraction, and hence as $\Kr(T)$
for some rational tangle $T$. Then $a^2+b^2$ (with the above remark on
signs) is the determinant of the achiral knot shown in \eqref{Tsum}.
\qed

Now we modify the second part of the proof to deduce theorem \ref{Th1}.
In the following we use Conway's notation for tangle sum and product.
(See for example again \cite[\S 2.3]{Adams} for a detailed description.)

\proof[of theorem \ref{Th1}] Let $n=p^2+q^2$. If $q=0$ then
$K=T(2,p)\#T(2,-p)$ ($T(2,p)$ denoting the $(2,p)$-torus knot)
is an easy example, so let $q\ne 0$. Krebes shows
that for any pair $(p,q)$ with at least one of $p$ and $q$ odd 
there is a(n arborescent) tangle $T$ with $\Kr(T)=
\frac{p}{q}$. In fact, $T$ can be chosen to be the connected sum of a
rational tangle and a knot of the type $T(2,p)$ (which can be done
in a way the tangle remains alternating).
Then again consider the knot in \eqref{Tsum}
(it is a knot because of the parities of $p$ and $q$), and from the
proof of theorem \ref{Th2} one sees
that it has the desired determinant $n$.

The knots constructed in \eqref{Tsum} then are all alternating.
It remains to show that they can be made prime (possibly sacrificing
alternation). If $(p,q)=1$, then $K$ is rational, and hence prime.
Thus let $n=(p,q)>1$. Then we can choose $T$ to be
\begin{eqn}\label{TT'}
T\quad=\quad\diag{7mm}{4.6}{2.9}{
  \pictranslate{-0.9 0}{
  \pictranslate{5 0.9}{
    \picrotate{90}{
      \lbraid{1 3.1}{2 2}
      \lbraid{0.35 d}{0.7 d}
      \lbraid{0.35 1.05}{0.7 d}
      \lbraid{0.35 1.75}{0.7 d}
      \pictranslate{-0.3 0}{
        \piccurve{1 2.1}{1 2.5}{1.5 2.5}{1.5 1.1}
        \piccurve{1 0}{1 -0.4}{1.5 -0.4}{1.5 1.1}
      }
      \piccurve{2 2.1}{1.9 1.6}{1.5 1}{2 -0.3}
    }
  }
  \picfilledbox{1.9 1.9}{1 1}{$T'$}
  \picputtext{3.95 0.8}{$
    \underbrace{\rx{2.1\unitlength}}_{n\scbox{ half-twists}}$}
  }
}\,,
\end{eqn}
that is, in Conway's notation $T=T'\cdot(0,n)$ with a tangle $T'$
being a rational tangle $a'/b'$ with $a'=p/n$, $b'=q/n$.

Now replace the $0$-tangle in \eqref{TT'} by the (flipped)
``KT-grabber'' tangle $KT$ in \cite{Bleiler,Flapan}. By the same
argument as in \cite{Flapan} (or see also \cite{KL,Van}), the tangle
$T_1=T'\cdot (KT\cdot 0,n)$ becomes prime, and hence so is then the
knot $K_1=\overline{T_1\cdot T_1}$ by proposition 1.3 of
\cite{Bleiler} (bar denotes tangle closure). As
in \cite{Bleiler}, $K_1$ and $K$ have the same Alexander
polynomial, so in particular the same determinant. \qed

% NEED TO REPL KT BY `KT\cdot 0' IF KT HAS CONN OF \infty RATHER
% THAN 0, CHECK \cite{Bleiler} !!!

\begin{exam}
To demonstrate the elegance of theorem \ref{Th1} as a chirality
criterion, we remark that among the prime knots of $\le 10$
crossings (denoted henceforth according to Rolfsen's tables
\cite[appendix]{Rolfsen}) there are 6 chiral knots with
self-conjugate HOMFLY polynomial~-- $9_{42}$, $10_{48}$,
$10_{71}$, $10_{91}$, $10_{104}$ and $10_{125}$, and this method shows
chirality of four of them~-- $9_{42}$, $10_{71}$, $10_{104}$ and
$10_{125}$, including the two examples ($9_{42}$ and $10_{71}$) where
additionally even the Kauffman polynomial is self-conjugate.
(For $9_{42}$ and $10_{125}$ the congruence modulo 4 is violated, so
that, as remarked on several other places, the signature works as well.)
\end{exam}

\begin{rem}
Since slice knots have square determinant,
it also follows that if there exists a rational knot
$S(p,q)$ which is at the same time achiral and slice, then it will
correspond to a Pythagorean triple, that is, be of the Schubert form
$S\bigl((m^2+n^2)^2,2mn(m^2-n^2)\bigr)$ with $m$ and $n$ coprime.
\end{rem}

With regard to theorem \reference{Th1}, we conjecture an
analogous statement to hold for links.

\begin{conj}\label{cj1}
An even natural number $n$ occurs as determinant of an
achiral link if and only if $n$ is the sum of two squares.
\end{conj}

As a remark on links, note that by the above description of numbers
which are sums of two squares, this set is closed under multiplication,
corresponding on the level of determinants of links to taking
connected sums. Thus it would suffice to prove conjecture \ref{cj1}
just for prime links.

More number theoretic results on the square representations (which
by the said above can also be transcribed knot-theoretically) may be
found in \cite{Kano,Williams}.

% \noindent{\bf Acknowledgements.} I would like to thank to
% Kenneth Williams for some helpful remarks.

% \section{Some generalizations and problems\label{par4}}

\subsection{Determinants of prime alternating achiral knots\label{Spa}}

Theorem \reference{Th1} naturally suggests the question in how
far the properties alternation and primeness can be combined
when realizing a sum of two squares as determinant of achiral knots.

In this case, the situation is much more difficult, though.
It is easy to see that not every (odd) sum of 2 squares can
be realized. The first (and trivial) example ist $1$, since
the only alternating knot with such determinant is the unknot
\cite{Crowell2}, and it is by definition not prime. However, there
are further examples.

\begin{propo}\label{prp1}
Let $n=9$ or $n=49$. Then
there is no prime alternating achiral knot of determinant $n$.
\end{propo}

By the above cited result of Crowell, one has a bound on the
crossing number of an alternating knot of given determinant, so
could check for any $n$ in finite time whether it is realized or not.
This renders the check for $n=9$ easy. However, the estimate we obtain
from Crowell's inequality is intractable in any practical sense for
$n=49$. We give an improvement of Crowell's result for achiral
alternating links, which, although not completely sharp, is enough
for our purpose. 

For the understanding and the proof of this result
we recall some standard terminology for knot diagrams.

\begin{defi}
The diagram on the left of figure \reference{figtan}
is called \em{connected sum} $A\# B$ of the diagrams $A$ and $B$.
If a diagram $D$ can be represented as the connected sum of 
diagrams $A$ and $B$, such that both $A$ and $B$ have at least one
crossing, then $D$ is called \em{disconnected} (or composite), else
it is called \em{connected} (or prime). Equivalently, a diagram
is prime if any closed curve intersecting it in exactly two points,
does not contain a crossing in one of its in- or exterior.

If a diagram $D$ can be written as $D_1\# D_2\#\,\cdots\,\# D_n$,
and all $D_i$ are prime, then they are called the \em{prime
(or connected) components/factors} of $D$.
\end{defi}

\begin{defi}
The diagram is \em{split}, if there is a closed curve not
intersecting it, but which contains parts of the diagram
in both its in- and exterior.
\end{defi}

\begin{figure}[htb]
\[
\diag{6mm}{3}{2}{
  \piccirclearc{1.8 1}{0.5}{-120 120}
  \picfilledcircle{1 1}{0.8}{$A$}
}\, \#\,
\diag{6mm}{3}{2}{
  \piccirclearc{1.2 1}{0.5}{60 300}
  \picfilledcircle{2 1}{0.8}{$B$}
}\quad =\quad
\diag{6mm}{4}{2}{
  \piccirclearc{2 0.5}{1.3}{45 135}
  \piccirclearc{2 1.5}{1.3}{-135 -45}
  \picfilledcircle{1 1}{0.8}{$A$}
  \picfilledcircle{3 1}{0.8}{$B$}
} 
\]
\caption{\label{figtan}}
\end{figure}

By \cite{Menasco} an alternating link is prime/split iff any
alternating diagram of it is so.

\begin{defi}
A crossing $q$ in a link diagram $D$ is called \em{nugatory}, if
there is a closed (smooth) plane curve $\gm$
intersecting $D$ transversely in $q$ and nowhere else.
A diagram is called \em{reduced} if it has no nugatory crossings.
\end{defi}

By \cite{Kauffman2,Murasugi3,Thistle}, each alternating reduced diagram
is of minimal crossing number (for the link it represents).

\begin{defi}
A {\em flype} is a move on a diagram shown in figure \reference{fig1}.

\begin{figure}[htb]
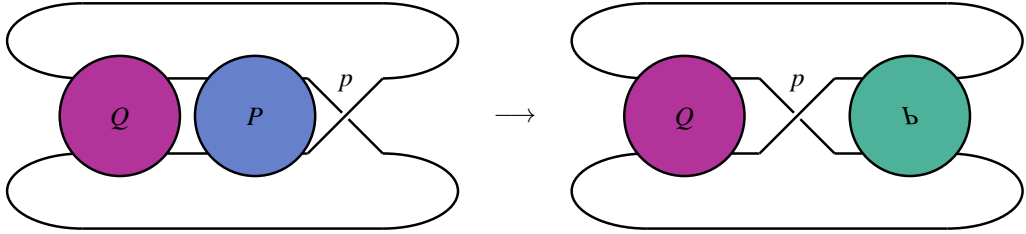

\[ %begin{equation}\label{flype}
\diag{1cm}{6}{3}{
    \pictranslate{3 1.5}{
       \picputtext[d]{1.5 0.4}{$p$}
       \picmultigraphics[S]{2}{1 -1}{
           \picmultiline{0.12 1 -1.0 0}{2 -0.5}{1 0.5}
           \picmultigraphics[S]{2}{-1 1}{
                \picellipsearc{-2 -1.0}{1 0.5}{90 270}
           }
           \picline{-2 -1.5}{2 -1.5}
           \picline{-2 -0.5}{1 -0.5}
      }
      \picfillgraycol{0.4 0.5 0.8}
      \picfilledcircle{0.3 0}{0.8}{$P$}
      \picfillgraycol{0.7 0.2 0.6}
      \picfilledcircle{-1.5 0}{0.8}{$Q$}
   }
}\quad\lra\quad
\diag{1cm}{6}{3}{
    \pictranslate{3 1.5}{
       \picputtext[d]{0 0.4}{$p$}
       \picmultigraphics[S]{2}{1 -1}{
           \picmultiline{0.12 1 -1.0 0}{0.5 -0.5}{-0.5 0.5}
           \picmultigraphics[S]{2}{-1 1}{
                \picellipsearc{-2 -1.0}{1 0.5}{90 270}
                \picline{-2 -0.5}{-0.5 -0.5}
           }
           \picline{-2 -1.5}{2 -1.5}
      }
      \picscale{1 -1}{
           \picfillgraycol{0.3 0.7 0.6}
           \picfilledcircle{1.5 0}{0.8}{$P$}
      }
      \picfillgraycol{0.7 0.2 0.6}
      \picfilledcircle{-1.5 0}{0.8}{$Q$}
   }
}
\] %end{equation}
\caption{A flype near the crossing $p$\label{fig_}}
\end{figure}
\end{defi}

By the fundamental work of \cite{MenThis}, for two alternating diagrams
of the same alternating link, there is a sequence of flypes
(and $S^2$-moves) taking the one diagram into the other.

\begin{defi}
To define the \em{sign} of a crossing in a link diagram,
choose an orientation of the link. The sign is then given as follows:
\begin{eqn}\label{bla}
\begin{array}{cc}
\diag{9mm}{1}{1}{
\picmultivecline{0.18 1 -1.0 0}{1 0}{0 1}
\picmultivecline{0.18 1 -1.0 0}{0 0}{1 1}
}
&
\diag{9mm}{1}{1}{
\picmultivecline{0.18 1 -1 0}{0 0}{1 1}
\picmultivecline{0.18 1 -1 0}{1 0}{0 1}
}
\\
\,+\,
&
\,-\,
\end{array}\,.
\end{eqn}
A crossing of sign $+$ we call \em{positive}, and a crossing
of sign $-$ we call \em{negative}. Note, that the definition requires
a link orientation, but for a \em{knot} it is independent on which of
its both possible choices is taken.

The \em{writhe} is the sum of the signs of all crossings in a diagram.
It is invariant under simultaneous reversal of orientation of
\em{all} components of the diagram, so is in particular well-defined
for unoriented \em{knot} diagrams. It may, however, change
if some (but not all) components of a link diagram are reverted.
\end{defi}

Crowell's result about the determinant of alternating links is

\begin{theorem}(Crowell)
If $L$ is a non-split alternating link of $n$ crossings, then
$\dt(L)\ge n$, and if $L$ is not the $(2,n)$-torus link, then
$\dt(L)\ge 2n-3$.
\end{theorem}

We will show

\begin{propo}\label{prp2}
If $L$ is an alternating non-split achiral link of $2n$ crossings,
then $\dt(L)\ge n(n-3)$.
\end{propo}

Since we use the checkerboard colorings for the proof, our
result holds for the most general notion of achirality for links~--
we allow the isotopy taking a link to its mirror image to interchange
components and/or preserve or reverse their orientations in an
arbitrary way. We will define, however, checkerboard colorings
later, in \S\reference{par5}, so that we defer the proof
of proposition \reference{prp2} to that later stage.

\begin{rem}
The condition the crossing number of $L$ to be even is no
restriction. The crossing number of any alternating achiral (in the
most general sense, as remarked after Proposition \reference{prp2})
link diagram is even by \cite{MenThis}, since flypes preserve the
writhe of the alternating diagram, and reversal of any single component
alters the writhe by a multiple of $4$ (any two components have an
even number of common crossings by the Jordan curve theorem).
Thus the writhe must be even, and hence so must be the crossing number.
\end{rem}

\proof[of proposition \reference{prp1}]
This is now feasible. Check all the alternating achiral knots in the 
tables of \cite{KnotScape} up to 16 crossings. \qed

There is some possibility that the values of proposition
\reference{prp1} are indeed the only exceptions.

\begin{conj}\label{Cjsq}
Let $n$ be an odd natural number. Then $n$ is the determinant of a
prime alternating achiral knot if and only if $n$
is the sum of two squares and $n\not\in\{1,\ 9,\ 49\}$.
\end{conj}

At least this is true up to $n\le 2000$. A full confirmation of
this conjecture is so far not possible, but we obtain a complete
statement for $n$ being a square.

\begin{theorem}\label{Thsq}
Let $n$ be an odd square. Then $n$ is the determinant of a prime
alternating achiral knot if and only if $n\not\in\{1,\ 9,\ 49\}$.
\end{theorem}

Recall, that a knot $K$ is called
\em{strongly achiral}, if it admits an embedding into $S^3$ pointwise
fixed by the (orientation-reversing) involution $(x,y,z)\mapsto
(-x,-y,-z)$. Again dependingly on the effect of this involution
on the orientation of the knot we distinguish between
strongly $+$achiral and strongly $-$achiral knots.

Let $\cL=\bQ[t,t^{-1}]$ be the Laurent polynomial ring in one
variable. For $F,G\in \cL$ write $F\sim G$ if $F$ and $G$ differ by
a multiplicative unit in $\bZ[t,t^{-1}]$, that is, $F(t)=\pm t^
nG(t)$ for some $n\in\bZ$.

The result of \cite{HarKaw} is the following.

\begin{theorem}(\cite{HarKaw})\label{thHK}
If $K$ is strongly negative amphicheiral, then $\Dl(t^2)\sim
F(t)F(-t)$ for some $F\in\cL$ with $F(-t)\sim F(t^{-1})$. 
If $K$ is strongly positive amphicheiral, then $\Dl(t)=F(t)^2$ for
some $F\in\cL$ with $F(t)\sim F(t^{-1})$.
\end{theorem}

This theorem will not be used here, but provides some heuristics
for the proof of theorem \reference{Thsq}, and will come to more
detailed mention later, so it is possibly appropriate to
introduce it here.

\proof[of theorem \reference{Thsq}]
Since $n$ is a square, it is suggestive by the result of
\cite{HarKaw} to consider strongly $+$achiral knots as
candidates to realize $n$ as determinant. We consider diagrams
of the type
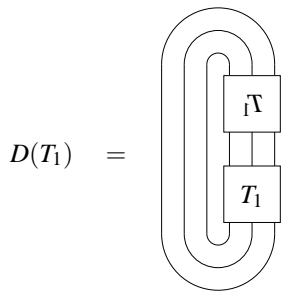
\[
% \begin{figure}[h]
% \fbox{
D(T_1)\quad=\quad
\diag{3mm}{5}{11}{
\pictranslate{3 2.0}{
   \picPSgraphics{0 setlinejoin}
   % \lbraid{0.5 0.5}{1 1}
   % \picline{2 0}{2 1}
   \picellipsearc{-0.5 0}{0.5 0.5}{180 0} 
   \picellipsearc{-0.5 0}{1.5 1.5}{180 0} 
   \picellipsearc{-0.5 0}{2.5 2.5}{180 0} 
   %\piccirclearc{2.5 0}{0.5}{180 0} 
   \picline{-1.0 0}{-1.0 7.5}
   \picline{-2.0 0}{-2.0 7.5}
   \picline{-3.0 0}{-3.0 7.5}
   \picline{ 1.0 0}{ 1.0 7.5}
   \picline{ 2.0 0}{ 2.0 7.5}
   \picline{ 0.0 0}{ 0.0 7.5}
   % \picline{3 0}{3 1}
%    \picmultigraphics{2}{0 3.5}{
%       \picline{0 1}{0 2}
%       \picline{2 2}{2 3}
%       \rbraid{1.5 1.5}{1 1}
%       \lbraid{0.5 2.5}{1 1}
%       % \lbraid{2.5 2.5}{1 1}
%    }
   \pictranslate{0 5.5}{
%       \picline{0 1}{0 2}
%       % \picline{3 1}{3 2}
%       \rbraid{1.5 1.5}{1 1}
      \picellipsearc{-0.5 2}{0.5 0.5}{0 180} 
      \picellipsearc{-0.5 2}{1.5 1.5}{0 180} 
      \picellipsearc{-0.5 2}{2.5 2.5}{0 180} 
%       %\piccirclearc{2.5 2}{0.5}{0 180} 
   }      
   \pictranslate{1 5.25}{\picscale{-1 1}{
   \picfilledbox{0 0.5}{2.5 2.5}{$T_1$}}
   \picfilledbox{0 -3.5}{2.5 2.5}{$T_1$}}
   % \picline{0 10}{0 11}
   % \picline{1 10}{1 11}
   % \lbraid{2.5 10.5}{1 1}
}}
% }
\]
Define a pairing $<\,T_1,T_2\,>$ on the diagram algebra
$DS_3(A=\sqrt{i})$ (see \cite{Kauffman}) by table \reference{tb3}
(compare also to the pairing $<\,,\,>_3$ in \S4 of \cite{det}).

\begin{table}
\[
\begin{array}{c*{5}{|c}}
\diag{3mm}{2}{2}{
  \picline{2 0}{0 2}
  \picputtext{0.3 d}{$T_1$}
  \picputtext{1.9 1.5}{$T_2$}
}             & \geni   & \genii  & \geniii & \geniv  & \genv   \\[3mm]
\hline
\rb{\geni}    &       0 &       0 &       0 &       1 &       1 \\[1mm]
\hline
\rb{\genii}   &       0 &       1 &       0 &       0 &       0 \\[1mm]
\hline
\rb{\geniii}  &       0 &       0 &       1 &       0 &       0 \\[1mm]
\hline
\rb{\geniv}   &       1 &       0 &       0 &       0 &       1 \\[1mm]
\hline
\rb{\genv}    &       1 &       0 &       0 &       1 &       0 \\[1mm]
\end{array}
\]
\caption{\label{tb3}}
\end{table}

Then $\dt(D(T_1))=<T_1,T_1>$. Let $T$ be a tangle
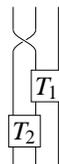
\[
\diag{3mm}{4}{7}{
   \picPSgraphics{0 setlinejoin}
   \picline{ 1.0 0}{ 1.0 5.0}
   \picline{ 2.0 0}{ 2.0 5.0}
   \picline{ 3.0 0}{ 3.0 7.0}
   \picline{ 1.0 6.0}{ 1.0 7.0}
   \picline{ 2.0 6.0}{ 2.0 7.0}
   \rbraid{ 1.5 5.5}{1.0 1.0}
   \pictranslate{2 0}{
   \picfilledbox{0.5 3.5}{1.4 1.4}{$T_1$}
   \picfilledbox{-0.5 1.5}{1.4 1.4}{$T_2$}
   %\picfilledbox{-0.5 5.5}{1.4 1.4}{$T_3$}
   }
}
\]
with $T_1=A\,\AT+B\,\BT$ and $T_2=X\,\AT+Y\,\BT$.
Then
\[
T=a\enspace\geni +b\enspace\genii +c\enspace\geniii + d\enspace\geniv+ e\enspace\genv
\]
with
\begin{eqnarray*}
a & = & XA \\
b\enspace =\enspace d & = & BX \\
c & = & AX+BY+AY \\
e & = & BY\,,
\end{eqnarray*}
and
\[
<\,T,T\,>\,\,  =\, \left[ (X+Y)(A+B) \right]^2\,.
\]
Whenever $(X,Y)$ and $(A,B)$ are relatively prime, and $X,Y,A,B>0$,
one can substitute rational tangles for $T_{1,2}$ obtaining
a prime alternating diagram of a strongly $+$achiral knot.
Setting $X=B=1$ and varying $Y$ and $A$, we see that we can cover
all cases when $n=p^2$ with $p$ composite.

Since we dealt with $p=1$, it remains to consider $p$ prime.
If $p\equiv 1\bmod 4$, then \eqref{*} shows that $n$ has a
non-trivial representation as sum of two squares, which then
must be coprime. In this case there is an achiral rational knot
of determinant $n$.

Thus assume $n=p^2$ with $p\equiv 3\bmod 4$ prime. We show now
that almost all (not necessarily prime) $p\equiv 3\bmod 4$ can be
realized.

Consider diagrams $D(T)$ for $T$ of the form
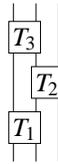
\[
\diag{3mm}{4}{7}{
   \picPSgraphics{0 setlinejoin}
   \picline{ 1.0 0}{ 1.0 7.0}
   \picline{ 2.0 0}{ 2.0 7.0}
   \picline{ 3.0 0}{ 3.0 7.0}
   \pictranslate{2 0}{
   \picfilledbox{0.5 3.5}{1.4 1.4}{$T_2$}
   \picfilledbox{-0.5 1.5}{1.4 1.4}{$T_1$}
   \picfilledbox{-0.5 5.5}{1.4 1.4}{$T_3$}
   }
}
\]
with $T_1=X\,\AT+Y\,\BT$, $T_2=A\,\AT+B\,\BT$ and $T_3=C\,\AT+D\,\BT$.
We find after multiplying out the polynomial and some manipulation
\[
<\,T,T\,>\,\, = \,\left[ X(DA+BC) + Y(BD+AC) \right]^2\,.
\]

Set $X=1$ and let $Y=k$ vary. The rest is done by choosing
small special values for $A,B,C,D$.

\def\Ta{\diag{3mm}{2}{2.3}{
  \piccurve{1.5 2}{1.1 1.7}{0.5 0.8}{0.5 d}
  \piccurveto{0.5 -0.2}{1.5 -0.2}{1.5 0.5}
  \piccurveto{1.5 0.8}{0.9 1.7}{0.5 2}
  \picline{0 0.5}{2 0.5}
} }
\def\Tb{\diag{3mm}{1}{2.3}{
  \piccurve{1 2}{-0.2 2}{-0.2 0}{1 0}
  \piccurve{0 2}{1.2 2}{1.2 0}{0 0}
} }
\def\Tc{\diag{3mm}{2}{4.3}{
  \pictranslate{1 0}{
    \picmultigraphics[S]{2}{-1 1}{
      \piccurve{0.5 4}{-1.7 3}{0.5 2.5}{0.5 1.9}
      \piccurveto{0.5 1.5}{-0.5 1}{-0.5 0.5}
      \piccurveto{-0.5 0.2}{-0.3 0}{0 0}
    }
  }
  \picline{0 0.5}{2 0.5}
} }
\def\Td{\diag{3mm}{3}{3.3}{
  \pictranslate{2 1.5}{
    \picmultigraphics[S]{2}{1 -1}{
      \piccurve{1 1.5}{-1 0.5}{1 0.1}{1 -0.6}
      \piccurveto{1 -0.8}{0.5 -1.5}{0 -1.5}
      \piccurveto{-0.9 -1.5}{-1.2 0.8}{-1.4 1.2}
    }
  }
} }
\def\Te{\diag{3mm}{1}{4.3}{
  \piccurve{1 0}{-1 1}{0.9 1.3}{0.9 2}
  \piccurveto{1 2.5}{-1 3}{1 4}
  \piccurve{0 0}{2 1}{0.1 1.3}{0.1 2}
  \piccurveto{0 2.5}{2 3}{0 4}
} }
\[
\begin{array}{*9{c|}c}
         & D & B & A & C & DA+BC & BD+AC & \sqrt{<\,T,T\,>} & T_2 & T_3 \\
\hline \ry{1.1em}
\xd{172} & 2 & 3 & 2 & 1 &   7   &   8   &   7+8k           & \Ta & \Tb \\
\xd{173} & 2 & 7 & 2 & 1 &  11   &  16   &  11+16k          & \Tc & \Tb \\
\xd{174} & 4 & 3 & 4 & 1 &  19   &  16   &  19+16k          & \Td & \Te 
\end{array}
\]

Examples of the knots thus obtained are given in figure
\reference{figg} (for simplicity, just the plane curves are drawn).

\begin{figure}[htb]
\[
\begin{array}{c@{\qquad}c@{\qquad}c}
\epsfsv{3cm}{t-ex1} & 
\epsfsv{3.1cm}{t-ex2} & 
\epsfsv{3cm}{t-ex3} \\
\ry{6mm} \xd{172},\ Y=4 & \xd{173},\ Y=3 & \xd{174},\ Y=3
\end{array}
\]
\caption{\label{figg}}
\end{figure}

All diagrams (and hence knots \cite{Menasco})
are prime for $k\ge 1$. Thus the only cases remaining
to check are for $\sqrt{n}=p\in\{3,7,11,19\}$. For $p=11$
we have $10_{123}$, and for $p=19$ we check the knots in the tables
of \cite{KnotScape}. We obtain the examples $12_{1019}$ (the
closure of the $5$-braid $(\sg_1\sg_2^{-1}
\sg_3\sg_4^{-1})^3$, with $\sg_i$ being the Artin
generators, as usual) and $14_{18362}$ (the closure of the $3$-braid
$\sg_1^2\sg_2^{-3}\sg_1^2\sg_2^{-2}\sg_1^3\sg_2^{-2}$). The cases
$p=3,\ 7$ were dealt with in proposition \reference{prp1}. \qed

\begin{rem}
We showed that in fact we can realize any $n$ stated in the
theorem by a strongly $+$achiral or rational knot, so
in particular by a strongly achiral knot, since an achiral rational
knot is known to be strongly $-$achiral (see \cite{HarKaw}).
It may be possible to exclude rational knots when allowing the
further exception $n=25$.
\end{rem}

To examine the general case of $n$ (not only perfect squares), one
needs to consider larger
series of examples. Because of \cite{HarKaw} the knots should
not (only) be strongly $+$achiral. A natural way to modify the
examples in the proof of theorem \reference{Thsq} to be $-$achiral
is to consider (braid type) closures of tangles like
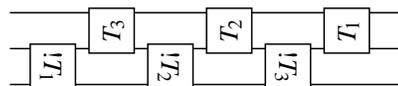
\[
\diag{6mm}{8}{2}{
  \pictranslate{3 1.2}{\picrotate{90}{
  \picline{-0.8 -4.5}{-0.8 4.2}
  \picline{0 -4.5}{0 4.2}
  \picline{0.8 -4.5}{0.8 4.2}
  \picfilledbox{0.4 1.95}{1.0 1.0}{$T_3$}
  \picfilledbox{0.4 -0.65}{1.0 1.0}{$T_2$}
  \picfilledbox{0.4 -3.25}{1.0 1.0}{$T_1$}
  \picrotate{180}{
    \picfilledbox{0.4 1.95}{1.0 1.0}{$!T_3$}
    \picfilledbox{0.4 -0.65}{1.0 1.0}{$!T_2$}
    \picfilledbox{0.4 -3.25}{1.0 1.0}{$!T_1$}
  }
  }}
}\,.
\]
We just briefly discuss this series to explain some of the occurring
difficulties.

Using the Kauffman bracket skein module coefficients
\begin{eqnarray*}
\diag{7mm}{1.5}{1}{
  \piclinewidth{49}
  \pictranslate{0.75 0.5}{\picrotate{90}{
    \picline{-0.4 -0.75}{-0.4 0.75}
    \picline{0.4 -0.75}{0.4 0.75}
    \picfilledbox{0 0}{1 1}{$T_1$}
  }}
} & = & A \A \, + \, B \B \,, \\
\diag{7mm}{1.5}{1}{
  \piclinewidth{49}
  \pictranslate{0.75 0.5}{\picrotate{90}{
    \picline{-0.4 -0.75}{-0.4 0.75}
    \picline{0.4 -0.75}{0.4 0.75}
    \picfilledbox{0 0}{1 1}{$T_2$}
  }}
} & = & C \A \, + \, D \B \,, \\
\diag{7mm}{1.5}{1}{
  \piclinewidth{49}
  \pictranslate{0.75 0.5}{\picrotate{90}{
    \picline{-0.4 -0.75}{-0.4 0.75}
    \picline{0.4 -0.75}{0.4 0.75}
    \picfilledbox{0 0}{1 1}{$T_3$}
  }}
} & = & X \A \, + \, Y \B \,, \\
\end{eqnarray*}
one finds an expression for $\dt(K)$ as polynomial in $A,B,C,D,X,Y$
as before, and after some manipulation arrives at $\dt(K)=f_1^2+f_2^2$
with 
\begin{eqn}\label{eqn2}
f_1(A,B,C,D,X,Y)\,=\,X(AD+BC)\,+\,Y(AC+BD)\qquad\mbox{and}\qquad
f_2(A,B,C,D,X,Y)\,=\,Y(AD-BC)\,.
\end{eqn}
(The correctness of the square decomposition is straightforward to
check, but for finding it it is helpful to notice that the substitutions
$Y=0$ and $A=C,\ B=D$ turn the knots into strongly $+$achiral ones,
which have square determinant.)

One can conclude from \eqref{eqn2} that no number of the form
$n=5p^2$ with $p\equiv 3\bmod 4$ prime can be written as
$f_1^2+f_2^2$ with $f_{1,2}$ as in \eqref{eqn2} for $A,B,C,D,X,Y\in\bN$,
unless $(A,B)$, $(C,D)$ or $(X,Y)$ is one of $(0,p)$, $(p,0)$ or
$(p,p)$. However, no alternating arborescent tangle has such
pair of Kauffman bracket skein module coefficients. Therefore, the
above series cannot realize these determinants. On the other hand,
the small cases in it (for $p\le 11$) are realized by knots
with Conway polyhedron $8^*$. 

Then one can consider more patters and write down more complicated
polynomials, each time having to show that each (at least sufficiently
large) number is realized by (at least some of) these polynomials.
Presently, such problems in number theory seem very difficult.
(One classic example is the determination of the numbers
$G(n)$ and $g(n)$ in Waring's problem, see for example \cite{DHL,%
Hooley,HardyWright}.) Therefore, conjecture \reference{Cjsq}
may be hard to approach as of now.

\subsection{Determinants of unknotting number one achiral knots\label{Su1}}

We conclude our results on sums of two squares by a related,
although somewhat auxiliary, consequence of the unknotting number
theorem of Lickorish \cite{Lickorish} and its refined version
given in \cite{unkn1}.

Let $u_{\pm}$ denote the \em{signed unknotting number}, the minimal
number of switches of crossings of a given sign to a crossing of the
reversed
sign needed to unknot a knot, or infinity if such an unknotting
procedure is not available (this is somewhat different from the
definition of \cite{CochLick}).

Thus a knot $K$ has $u_+(K)=1$ (resp.\ $u_-(K)=1$) if it can be
unknotted by switching a positive (resp.\ negative) crossing in some of
its diagrams.
% Furthermore, let $\lm$ be the linking form on $H_1(D_K)$. the
% homology group (over $\bZ$) of the double branched cover of $S^3$ over
% $K$ (whose order, as discussed in the introduction, is given by
% $\big|\Dl_K(-1)\big|$).

\begin{prop}\label{pp1}
Let $K$ be a knot with $u_+=u_-=1$ (for example, an achiral unknotting
number one knot). Then $\det(K)$ is the sum of two squares of coprime
numbers.
\end{prop}

\proof Clearly $\sg(K)=0$, so that any of the relevant
crossing changes does not alter the signature, and then by the
refinement of Lickorish's theorem given in \cite{unkn1}, we have
$\lm(g_\pm,g_\pm)=\pm2/\det(K)\in\bQ/\bZ$ for some generators
$g_\pm$ of $H_1(D_K)$. Here $\lm$ is the linking form on $H_1(D_K)$
(whose order, as discussed in the introduction, is given by
$|\Dl_K(-1)|$\,). Thus $2l^2=-2h^2$ for some $l,h\in\bZ_{\det(K)}^*$.
Then this group possesses square roots of $-1$.

The structure of the group $\bZ_n^*$ of units in $\bZ_n=\bZ/n\bZ$
is known; see e.g. \cite[exercise 1, \S 5, p.\ 41]{Zagier}.
{}From this structure and \eqref{r02}, we see that the number of
square roots of $-1$ in $\bZ_n^*$ (for $n>1$ odd) is identical
to $r_2^0(n)$. (D.\ Zagier remarked to me that one can in
fact give a natural bijection between the square roots of $-1$
in $\bZ_n^*$ and representations of $n$ as sum of coprime
squares.) \qed

% For $-1$ this number is closely related to the above mentioned
% quantity $r_2^0(n)$.
% and the result
% follows from the relation of the number of such roots to $r_2^0$
% remarked above. \qed

% The conditions may be somewhat technical, but they are often satisfied.
% For example the first condition is satisfied for all rational knots,
% and the second one for all achiral unknotting number one knots.

There are also several questions opened by proposition \reference{pp1}.
Having the inclusions
\[
\{\,\dt(K)\,:\,K\mbox{ achiral, }u(K)=1\,\}\,\subset\,
\{\,\dt(K)\,:\,u_+(K)=u_-(K)=1\,\}\,\subset\,
\{\,a^2+b^2\,:\,(a,b)=1,\enspace 2\nmid a+b\,\}
\]
(to introduce some notation, let $S_a$, $S_{\pm}$ and $S$ denote these
three sets from left to right), the first question is whether and/or
which one of these inclusions is proper (or not).

This seems much more difficult to decide than the proof of theorems
\reference{Th1} and \reference{Th2}. There is no such straightforward
procedure available to exhaust all values in $S$, and to show a proper
inclusion one will face the major problem of deciding about unknotting
number one. 

After a computer experiment with the knot tables and tools
available to me, the smallest $x\in S$ I could not decide to belong
to $S_\pm$ is 349, whereas the smallest possible $x\in S$ with
$x\not\in S_a$ is only 17. Contrarily to theorems \reference{Th1} and
\reference{Th2}, very many entries have been completed only by
non-alternating knots (which have smaller determinant than the
alternating ones of the same crossing number), and in fact we can
use the number 17 to show that for this problem non-alternating
knots definitely need to be considered.

\begin{exam}
We already quoted Crowell's result,
that for a given crossing number $n$, the $n$ crossing twist knot
has the smallest determinant $2n-3$ among alternating knots $K$ of
crossing number $n$, except if $n$ is odd, in which case the
only knot of smaller determinant is the $(2,n)$-torus
knot. (A more modern proof can be given for example by the Kauffman
bracket, similarly to proposition \reference{prp2}.)
Thus, except for the $(2,17)$-torus knot, any alternating
knot of determinant 17 has $\le 10$ crossings. A direct check shows that
the only such knots of $\sg=0$ are $8_3$ and $10_1$. However, $u(8_3)=2$
as shown by Kanenobu-Murakami \cite{KanMur}, and that $10_1$ cannot
simultaneously have $u_+=u_-=1$ follows by refining their method
(see \cite{AskS}). Thus there is no alternating knot of determinant 17
with $u_+=u_-=1$, and the inclusion
\[
S_{a\pm}\,:=\,\{\,\dt(K)\,:\,K\mbox{ alternating, }
u_+(K)=u_-(K)=1\,\}\,\subset\,S
\]
is proper. (Contrarily, there is a simple non-alternating knot,
$9_{44}$, with $u_+=u_-=1$ and determinant 17.)
Is it infinitely proper, i.e., is $|S\sm S_{a\pm}|=\infty\,$?
\end{exam}

It is suggestive that considering the even more restricted class
of rational knots, the inclusions are infinitely proper. We confirm this
for achiral unknotting number one rational knots.

\begin{prop} $\ds\big|\,\,S\sm
\{\,\dt(K)\,:\,K\mbox{ achiral and rational, }
u(K)=1\,\}\,\big|\,=\,\infty\,$. (In fact, this set
contains infinitely many primes.)
\end{prop}

\proof These knots were classified in \cite[corollary 2.2]{rational}
to be those with Conway notation $(n11n)$ and $(3(12)^n1^4(21)^n3)$.
It is easy to see that therefore the inclusion
\[
\{\,\dt(K)\,:\,K\mbox{ achiral and rational, }
u(K)=1\,\}\,\subset\,S
\]
is infinitely proper. For example, the determinant of both series grows
quadratically resp.\ exponentially in $n$, so that
\[
\sum_{\scbox{$\begin{array}{c}K=S(p,q)\\\mbox{achiral, }u(K)=1
\end{array}$}} \frac 1p\,<\,\infty\,,
\]
while $\ds\sum_{p\in S}\frac 1p\,=\,\infty$. (By Dirichlet already
$\ds\sum_{\scbox{$\begin{array}{c}p\equiv 1\,(4)\\\mbox{prime}
\end{array}$}}\frac 1p\,=\,\infty$\,, see \cite[Korollar, p.~46]
{Zagier}.) \qed

It appears straightforward to push the method of \cite{rational}
further to show the same also for rational knots with $u_+=u_-=1$
(although I have not carried out a proof in detail). By refining
Kanenobu-Murakami, we first show that there is a crossing of the
same sign unknotting the alternating diagram of the rational knot (see
\cite{AskS}). Then consider alternating rational knot diagrams with two
unknotting crossings (it does not even seem necessary to have them
any more of different sign), and apply the same argument as above,
using \cite[corollary 2.3]{rational} and the remark after its proof.
% (The statement in this remark is only suggested,
% and is yet to be proved.)

\section{Enumeration of rational knots by determinant\label{Senum}}

The results of \S\reference{S2} can be used to
enumerate rational knots by determinant. We have for example:

\begin{prop}\label{xp}
The number $c_n$ of achiral rational knots of
given determinant $n$ is given by
\[
c_n\,=\,\left\{\begin{array}{ll}
\myfrac 12\,{r_2^0(n)} & \mbox{ if $n>2$ odd,} \\
0 & \mbox{ else\,.}
\end{array}\right.
\]
\end{prop}

\proof Use
the fact that there is a bijective correspondence between the rational
tangle $T$ in a diagram \eqref{Tsum} of an achiral rational knot
$K$ and its Krebes invariant $p/q$ (with $p\ge q$ and $(p,q)=1$)
giving $\dt(K)=p^2+q^2$. \qed

% As such representations are clearly in one-to-one correspondence with
% the achiral rational knots, we obtain
{}From \eqref{r02} we obtain then

\begin{corr}
The number of achiral rational knots of given determinant $n$
is either zero or a power of two. \qed
\end{corr}

%\begin{rem} (after proof of theo 2)
As a practical application of the argument in the argument proving
proposition \reference{xp} we can consider the achiral rational knots
$(1 \dots 1)$ and $(3 1 \dots 1 3)$ (with the number list of even
length) and the tangles $T$ obtained from the halves of the palindromic
sequence. This way one arrives to a knot theoretical explanation of
the identities
\begin{eqn}
F_{2n+1}=F_n^2+F_{n+1}^2\quad\mbox{and}\quad
L_{2n+1}+2L_{2n}=L_n^2+L_{n+1}^2\,,
\end{eqn}
where $F_n$ is the $n$-th Fibonacci number ($F_1=1$, $F_2=1$,
$F_n=F_{n-1}+F_{n-2}$) and $L_n$ is the $n$-th Lucas number
($L_1=2$, $L_2=1$, $L_n=L_{n-1}+L_{n-2}$). Thus we have
% In particular these representations of odd index Fibonacci and Lucas
% numbers show

\begin{prop}
There are achiral rational knots with determinant $F_{2n+1}$ and
$L_{2n+1}+2L_{2n}$ for any $n$, or equivalently, any (prime)
number $4x+3$ does not divide $F_{2n+1}$ and $L_{2n+1}+2L_{2n}$. \qed
\end{prop} 
%of the above identities%

For $F_{2n+1}$ this is a task I remember from an old issue of
the Bulgarian journal ``Matematika''. Recently I found (by electronic
search) that it was conjectured in \cite{Thoro} and proved in
\cite{Yalavigi}.
%
% \begin{corr}
% Any $F_{2n+1}$ does not have any divisor of the form $4x+3$.
% \end{corr}
% 
% \proof Clearly it suffices to prove the assertion when $q=4x+3$ is
% prime. We proceed by induction on $n$. The case $n=0$ is trivial, and
% the identity shows that for any $n$, $F_{2n+1}\equiv x^2+y^2\bmod q$
% with one of $x^2$ and $y^2$ by induction assumption being non-zero.
% Therefore, we are done unless the other square is non-zero (modulo $q$)
% either. But then $x^2,\ y^2\in\bZ_q^*$, which does not have any element
% of order $4$, and hence $-x^2\in\bZ_q^*$ cannot be a square itself. \qed

A similar enumeration can be done for arbitrary rational knots
of given (odd) determinant $n$, and one obtains

\begin{prop}
The number of rational knots of determinant $n$ ($n>1$ odd),
counting chiral pairs once, is
\begin{eqn}\label{rr}
\frac{1}{4}\left\{\,\phi(n)+r_2^0(n)+2^{\om(n)}\right\}\,,
\end{eqn}
with $r_2^0(n)$ being as in \eqref{r20}, $\om(n)$ denoting the number
of different prime divisors of $n$ and $\phi(n)$ being Euler's totient
function.
\end{prop}

\proof We apply Burnside's lemma on the action of $\bZ_2\times
\bZ_2$ on $\bZ_n^*$ given by additive inversion in the first
component and multiplicative inversion in the second one. In \eqref{rr},
the second and third term in the braced expression come from counting
the square roots of $\mp 1$ in $\bZ_n^*$. These numbers
follow from the structure of this group $\bZ_n^*$,
as remarked in the proof of proposition \reference{pp1}. \qed

\begin{rem}
The functions $\om(n)$ and $\phi(n)$ are hard to calculate for
sufficiently large numbers $n$ by virtue of requiring the prime
factorization of $n$, but the expression in terms of these
classical number theoretical functions should be at least of theoretical
interest.
\end{rem}

Counting chiral pairs twice one has the somewhat simpler expression
\[
\frac{1}{2}\left\{\,\phi(n)+2^{\om(n)}\right\}\,.
\]

In a similar way one could attempt the enumeration by $c_p$ of
unknotting number one rational knots of determinant $p$ using
\cite{KanMur}, seeking again an expression in terms of classical
number theoretical functions. Obviously from the result of
Kanenobu-Murakami we have
\[
c_p\,\le\,2^{\om((p+1)/2)-1}+2^{\om((p-1)/2)-1}-1\,,
\]
with the powers of two counting the representations of $(p\pm1)/2$
as the product of two coprime numbers $n_{\pm}$ and $m_{\pm}$ up
to interchange of factors and the final `$-1$' accounting for the double
representation of the twist knot for $m_+=m_-=1$. However, to
obtain an exact formula, one encounters the problem that, beside
the twist knot, some other knot may arise from different
representations (although this does not occur often and the inequality
above is very often sharp). For example, for $p=985$ the knot
$S(985,288)=S(985,697)$ occurs for the representations $m_+=29$
and $m_-=12$. D.\ Zagier informed me that he has obtained a complete
description
of the duplications of the Kanenobu-Murakami forms when considering
$q$ in $S(p,q)$ only up to \em{additive} inversion in $\bZ_p^*$.
According to him, however, considering the (more relevant)
\em{multiplicative} inversion renders the picture too complicated and
number theoretically unilluminating.
% (In this case the duplication occurs among representations
% of different sign choice, but it is not \em{a priori} clear whether
% duplications may not occur even among representations of the same sign.)

\section{Spanning trees in planar graphs and checkerboard colorings%
\label{par5}}

Here we discuss an interpretation of our results of \S\reference{S3}
in graph theoretic terms.

\begin{theorem}\label{th5.1}
Let $n$ be an odd natural number. Then $n$ is the number of
spanning trees in a planar self-dual graph if and only if
$n$ is the sum of two squares.
\end{theorem}

The proof of this theorem relies on the following construction
linking graph and knot theory (see e.g. \cite{Kauffman}).
Given an alternating knot (or link) diagram $D$, we can associate
to it its checkerboard graph.

The \em{checkerboard coloring} of a link diagram is a map
\[
\{\mbox{\ \ regions of $D$\ \ }\}\,\to\,
\{\mbox{\ \ black, white\ \ }\}
\]
s.t. regions sharing an edge are always mapped to different colors.
(A \em{region} is called a connected component of the complement of the
plane curve of $D$, and an \em{edge} a part of the plane curve of the
diagram between two crossings.)

The \em{checkerboard graph} of $D$ is defined to have vertices
corresponding to black regions in the checkerboard coloring of $D$,
and an edge for each crossing $p$ of $D$ connecting the two black
regions opposite at crossing $p$ (so multiple edges between two
vertices are allowed).

This construction defines a bijection
\[
\{\mbox{\ \ alternating diagrams up to mirroring\ \ }\}\quad\llra\quad
\{\mbox{\ \ planar graphs up to duality\ \ }\}\,.
\]

Duality of the planar graph corresponds to switching
colors in the checkerboard coloring and has the effect of
mirroring the alternating diagram if we fix the sign
of the crossings so that each crossing looks like
$
\diag{4mm}{2}{2}{
  \picgraycol{1}
  \picfillgraycol{0.75}
  \picfill{
    \picline{0 0}{0 2}\piclineto{2 0}\piclineto{2 2}
  }
  \picgraycol{0}
  \picline{2 0}{1.1 0.9}
  \picline{2 0 x}{1.1 0.9 x}
  \picline{0 0}{2 2}
}
$ rather than
$
\diag{4mm}{2}{2}{
  \picgraycol{1}
  \picfillgraycol{0.75}
  \picfill{
    \picline{0 0}{0 2}\piclineto{2 0}\piclineto{2 2}
  }
  \picgraycol{0}
  \picline{2 0}{0 2}
  \picline{0 0}{0.9 d}
  \picline{1.1 d}{2 2}
}
$. Then we have

\begin{lemma}
$\dt(D)$ is the number of spanning trees in a checkerboard graph
of $D$ for any alternating link diagram $D$.
\end{lemma}

\proof By the Kauffman bracket definition of the Jones polynomial $V$,
for an alternating diagram $D$, the determinant $\dt(D)=|\Dl_D(-1)|=
|V_D(-1)|$ can be calculated as follows (see \cite{Krebes}).

Consider $\hat D\subset\bR^2$, the (image of) the associated
immersed plane curve(s). For each crossing (self-intersection) of
$\hat D$ there are 2 ways to splice it:
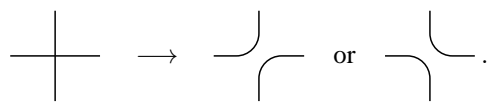
\[
\diag{6mm}{2}{2}{
  \picline{1 0}{1 1}
  \picline{2 1}{1 1}
  \picline{1 2}{1 1}
  \picline{0 1}{1 1}
} \quad\lra\quad
\diag{6mm}{2}{2}{
  \picarcangle{1 0}{1 1}{2 1}{0.5}
  \picarcangle{1 2}{1 1}{0 1}{0.5}
}\mbox{\quad or\quad}
\diag{6mm}{2}{2}{
  \picarcangle{1 2}{1 1}{2 1}{0.5}
  \picarcangle{1 0}{1 1}{0 1}{0.5}
}\,.
\]
We call a choice of splicing for each crossing a \em{state}.
Then $\dt(D)$ is equal to the number of states
so that the resulting collection of disjoint circles has only one
component (a single circle). We call such states \em{monocyclic}.

Let $\Gm$ be a spanning tree of the checkerboard graph $G$ of $D$.
Define a state $S(\Gm)$ as follows: for any edge $v$ in $G$ set
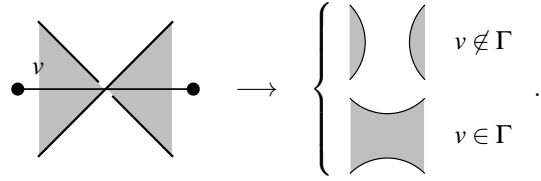
\[
\diag{9mm}{2}{2}{
  \picgraycol{1}
  \picfillgraycol{0.75}
  \picfill{
    \picline{0 0}{0 2}\piclineto{2 0}\piclineto{2 2}
  }
  \picgraycol{0}
  \picline{2 0}{1.1 0.9}
  \picline{2 0 x}{1.1 0.9 x}
  \picline{0 0}{2 2}
  \picfillgraycol{0}
  \picfilledcircle{2.3 1}{0.08}{}
  \picfilledcircle{-0.3 1}{0.08}{}
  \piclinewidth{40}
  \picline{-0.3 1}{2.3 1}
  \picputtext{0 1.3}{$v$}
}\qquad\lra\quad
\left\{
\begin{array}{cc}
\diag{5mm}{2}{2}{
  \picgraycol{1}
  \picfillgraycol{0.75}
  \picfill{
    \piccirclearc{-1 1}{1.41}{-45 45}
    \piccirclearc{3 1}{1.41}{135 -135}
  } 
  \picgraycol{0}
  \piccirclearc{-1 1}{1.41}{-45 45}
  \piccirclearc{3 1}{1.41}{135 -135}
} & v\not\in\Gm \\[5mm]
\diag{5mm}{2}{2}{
  \picgraycol{1}
  \picfillgraycol{0.75}
  \picfill{\picbox{1 1}{2 2}{}}
  \picfillgraycol{1}
  \picfill{
    \picline{0 2}{0 0}
    \piccirclearc{1 -1}{1.41}{45 135}
    \picline{2 0}{2 2}
    \piccirclearc{1 3}{1.41}{-135 -45}
  }
  \picgraycol{0}
  \piccirclearc{1 -1}{1.41}{45 135}
  \piccirclearc{1 3}{1.41}{-135 -45}
} & v\in\Gm \\
\end{array}
\right.\,.
\]
Then $S$ gives a bijection between monocyclic states of $D$ and
spanning trees of $G$. \qed

Since $\dt(K)=|\Dl_K(-1)|$ and it is known that for an $n$-component
link $K$, $(t^{1/2}-t^{-1/2})^{n-1}\mid \Dl_K(t)$, we have that
$2^{n-1}\mid \dt(K)$. Thus $\dt(K)$ is odd only if $K$ is a knot.
The converse is also true, since for a knot $K$ we have $\Dl_K(1)=1$,
and $\Dl_K(-1)\equiv \Dl_K(1)\bmod 2$.

\proof[of theorem \reference{th5.1}]
If $n$ is the number of spanning trees of a planar self-dual graph $G$,
then its associated alternating diagram $D$ is isotopic by $S^2$-moves
to its mirror image. Since $n$ is assumed odd, $D$ is a knot diagram.
Thus the number of spanning trees of $G$, which by the lemma is equal
to $\dt(D)$, is of the form $a^2+b^2$ by theorem \reference{Th1}.

Contrarily, assume that $n=a^2+b^2$. Take the checkerboard graph $G$
of the diagram in \eqref{Tsum} constructed in the proof of
theorem \ref{Th2}. This diagram has the property of being
isotopic to its mirror image by $S^2$-moves only (and no flypes),
so that its (self-dual) checkerboard graph $G$ is the one we sought.
\qed

Using checkerboard colorings, we will now give a proof
of proposition \reference{prp2}. We introduce first some more
standard notations.

\begin{defi}
Given a diagram $D$ and a closed curve $\gm$ intersecting $D$
in exactly four points, $\gm$ defines a \em{tangle decomposition}
of $D$.
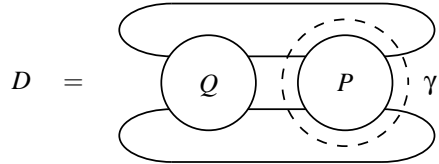
\[
D\quad=\quad
\diag{7mm}{6}{3}{
  \pictranslate{3 1.5}{
       \picmultigraphics[S]{2}{1 -1}{
           \picmultiline{0.12 1 -1.0 0}{2 -0.5}{1 0.5}
           \picmultigraphics[S]{2}{-1 1}{
                \picellipsearc{-2 -1.0}{1 0.5}{90 270}
           }
           \picline{-2 -1.5}{2 -1.5}
           \picline{-2 -0.5}{1 -0.5}
      }
  }
  \pictranslate{1.7 1.5}{
      \picfilledcircle{0 0}{0.9}{$Q$}
  }
  \picfilledcircle{4.3 1.5}{0.9}{$P$}
  \piclinedash{0.2}{0.5}
  \piccircle{4.3 1.5}{1.2}{}
  \picputtext[l]{5.8 1.5}{$\gm$}
}
\]
A \em{mutation} of $D$ is obtained by removing one of
the tangles in some tangle decomposition of $D$ and replacing it
by a rotated version of it by $180^\circ$ along the axis vertical to
the projection plane, or horizontal or vertical in the projection plane.
For example:
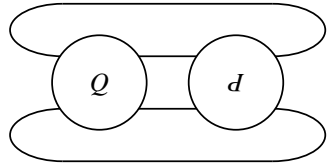
\[
\diag{7mm}{6}{3}{
  \pictranslate{3 1.5}{
       \picmultigraphics[S]{2}{1 -1}{
           \picmultiline{0.12 1 -1.0 0}{2 -0.5}{1 0.5}
           \picmultigraphics[S]{2}{-1 1}{
                \picellipsearc{-2 -1.0}{1 0.5}{90 270}
           }
           \picline{-2 -1.5}{2 -1.5}
           \picline{-2 -0.5}{1 -0.5}
      }
  }
  \pictranslate{1.7 1.5}{
    \picrotate{0}{
      \picfilledcircle{0 0}{0.9}{$Q$}
    }
  }
  \pictranslate{4.3 1.5}{
    \picrotate{180}{
      \picfilledcircle{0 0}{0.9}{$P$}
    }
  }
}
\]
(To make the orientations compatible, eventually the orientation
of either $P$ or $Q$ must be altered.)
$\gm$ is called the \em{Conway circle} for this mutation.
\end{defi}

Note that a flype can be realized as a sequence of mutations.

There is an evident bijection between the crossings of a diagram 
before and after applying a mutation, so that we can trace a crossing
in a sequence of mutations and identify it with its image in the
transformed diagram when convenient. In particular, we can do so
for a sequence of flypes.

\proof[of proposition \reference{prp2}]
We proceed by induction on $n$. For $n\le 1$ the claim is trivial.

Assume now $L$ have $2n$ crossings and be alternating and achiral,
and $n>1$. By \cite{MenThis}, there is a sequence of flypes
(and $S^2$-moves) taking an alternating diagram $D$ of $L$ into
its mirror image $!D$.

Fix a crossing $p$ in $D$ and let $p'$ be the crossing in $D$
whose trace under the flypes taking $D$ to $!D$ takes it
to the mirror image of $p$ in $!D$.

Since the only diagram in which both splicings of a crossing
give a nugatory crossing is the Hopf link diagram, for each
$p$ and $p'$ there is a splicing not producing a nugatory
crossing. We call such splicing a \em{non-nugatory splicing},
otherwise call the splicing \em{nugatory}.

We can distinguish the two splicings at $p$ according to
the colors of the regions in the checkerboard coloring they
join. We thus call the splicings \em{black} or \em{white}.

Since whether the black or white splicing gives a
nugatory crossing is invariant under flypes,
the choice of non-nugatory splicing at $p$ and $p'$ between
black or white splicing can be made to be the opposite.
Call $D_p$ (resp.\ $D_{p'}$) the diagrams obtained from $D$
after the so chosen splicing at $p$ (resp.\ $p'$), and $D'$
the diagram resulting after performing both splicings. Since
$D_p$ and $D_{p'}$ have no nugatory crossings, $D'$ is non-split.

We claim that $D'$ has no nugatory crossings. To see this, use that
$p$ and $p'$ join two pairs of regions of opposite color
in the checkerboard coloring of $D$. If $D'$ has a nugatory
crossing $q$, then there would be a closed plane curve $\gm$
intersecting $D'$ (transversely) only in $q$, and lying in some
(without loss of generality) white region of the
checkerboard coloring of $D'$.
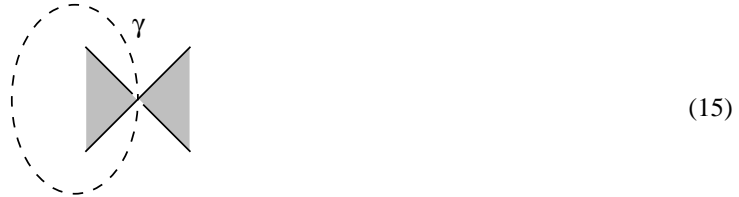
\begin{eqn}\label{gmeq}
\diag{7mm}{3}{4}{
  \pictranslate{1 1.2}{
    \picgraycol{1}
    \picfillgraycol{0.75}
    \picfill{
      \picline{0 0}{0 2}\piclineto{2 0}\piclineto{2 2}
    }
    \picgraycol{0}
    \picline{2 0}{1.1 0.9}
    \picline{2 0 x}{1.1 0.9 x}
    \picline{0 0}{2 2}
  }
  \piclinedash{0.2}{0.5}
  \picellipse{0.8 2.2}{1.2 1.8}{}
  \picputtext{2.0 3.6}{$\gm$}
}
\end{eqn}
But then $\gm$ would persist by undoing the splicing joining
the two black regions, and thus $q$ would be nugatory in one of
$D_p$ or $D_{p'}$, too, a contradiction.

Since we need the plane curve argument later, let us for convenience
call a curve $\gm$ through black (resp.\ white) regions
of $D$ intersecting $D$ in crossings $c_1,\dots,c_n$
a \em{black (resp.\ white) curve} through $c_1,\dots,c_n$.

We also claim that $D'$ depicts the mutant of an achiral link. This
follows, because the flypes carrying $D$ to $!D$ also carry $D'$ to
$!D'$, modulo mutation.
To see this, the only problematic case to consider is when a flype at
$p$ must be performed. Then
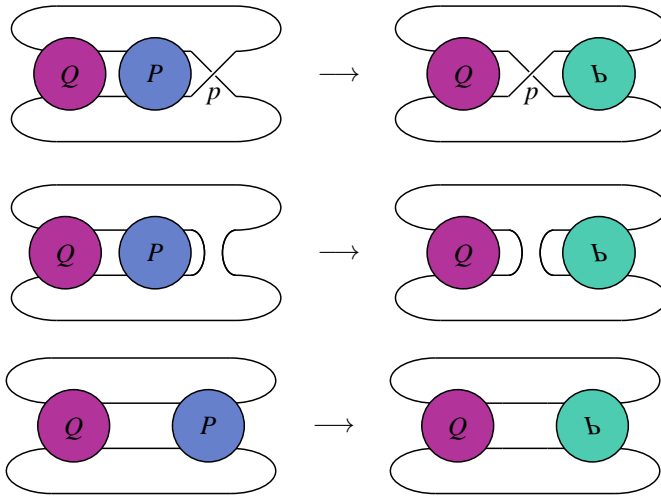
\[
\diag{6mm}{6}{3}{
    \pictranslate{3 1.5}{
      \picmultigraphics[S]{2}{1 -1}{
           \picmultiline{0.12 1 -1.0 0}{2 -0.5}{1 0.5}
           \picmultigraphics[S]{2}{-1 1}{
                \picellipsearc{-2 -1.0}{1 0.5}{90 -90}
           }
           \picline{-2 -1.5}{2 -1.5}
           \picline{-2 -0.5}{1 -0.5}
      }
      \picputtext[u]{1.5 -0.3}{$p$}
      \picfillgraycol{0.4 0.5 0.8}
      \picfilledcircle{0.2 0}{0.8}{$P$}
      \picfillgraycol{0.7 0.2 0.6}
      \picfilledcircle{-1.7 0}{0.8}{$Q$}
   }
}\quad\lra\quad
\diag{6mm}{6}{3}{
    \pictranslate{3 1.5}{
       \picmultigraphics[S]{2}{1 -1}{
           \picmultiline{0.12 1 -1.0 0}{0.5 -0.5}{-0.5 0.5}
           \picmultigraphics[S]{2}{-1 1}{
                \picellipsearc{-2 -1.0}{1 0.5}{90 -90}
                \picline{-2 -0.5}{-0.5 -0.5}
           }
           \picline{-2 -1.5}{2 -1.5}
      }
      \picputtext[u]{0 -0.3}{$p$}
      \picscale{1 -1}{
           \picfillgraycol{0.3 0.8 0.7}
           \picfilledcircle{1.5 0}{0.8}{$P$}
      }
      \picfillgraycol{0.7 0.2 0.6}
      \picfilledcircle{-1.5 0}{0.8}{$Q$}
   }
}
\]
turns into
\[
\diag{6mm}{6}{3}{
    \pictranslate{3 1.5}{
       \picmultigraphics[S]{2}{1 -1}{
	   \picellipsearc{2 0}{0.3 0.5}{90 270}
	   \picellipsearc{1 0}{0.3 0.5}{270 90}
           % \picmultiline{0.12 1 -1.0 0}{2 -0.5}{1 0.5}
           \picmultigraphics[S]{2}{-1 1}{
                \picellipsearc{-2 -1.0}{1 0.5}{90 -90}
           }
           \picline{-2 -1.5}{2 -1.5}
           \picline{-2 -0.5}{1 -0.5}
      }
      \picfillgraycol{0.4 0.5 0.8}
      \picfilledcircle{0.2 0}{0.8}{$P$}
      \picfillgraycol{0.7 0.2 0.6}
      \picfilledcircle{-1.8 0}{0.8}{$Q$}
   }
}\quad\lra\quad
\diag{6mm}{6}{3}{
    \pictranslate{3 1.5}{
       \picmultigraphics[S]{2}{1 -1}{
	   \picellipsearc{0.5 0}{0.3 0.5}{90 270}
	   \picellipsearc{-0.5 0}{0.3 0.5}{270 90}
           % \picmultiline{0.12 1 -1.0 0}{0.5 -0.5}{-0.5 0.5}
           \picmultigraphics[S]{2}{-1 1}{
                \picellipsearc{-2 -1.0}{1 0.5}{90 -90}
                \picline{-2 -0.5}{-0.5 -0.5}
           }
           \picline{-2 -1.5}{2 -1.5}
      }
      \picscale{1 -1}{
           \picfillgraycol{0.3 0.8 0.7}
           \picfilledcircle{1.5 0}{0.8}{$P$}
      }
      \picfillgraycol{0.7 0.2 0.6}
      \picfilledcircle{-1.5 0}{0.8}{$Q$}
   }
}
\]
or
\[
\diag{6mm}{6}{3}{
    \pictranslate{3 1.5}{
       \picmultigraphics[S]{2}{1 -1}{
           \picmultiline{0.12 1 -1.0 0}{2 -0.5}{1 0.5}
           \picmultigraphics[S]{2}{-1 1}{
                \picellipsearc{-2 -1.0}{1 0.5}{90 -90}
           }
           \picline{-2 -1.5}{2 -1.5}
           \picline{-2 -0.5}{1 -0.5}
      }
      \picfillgraycol{0.4 0.5 0.8}
      \picfilledcircle{1.5 0}{0.8}{$P$}
      \picfillgraycol{0.7 0.2 0.6}
      \picfilledcircle{-1.5 0}{0.8}{$Q$}
   }
}\quad\lra\quad
\diag{6mm}{6}{3}{
    \pictranslate{3 1.5}{
       \picmultigraphics[S]{2}{1 -1}{
           \picline{0.5 0.5}{-0.5 0.5}
           \picmultigraphics[S]{2}{-1 1}{
                \picellipsearc{-2 -1.0}{1 0.5}{90 -90}
                \picline{-2 -0.5}{-0.5 -0.5}
           }
           \picline{-2 -1.5}{2 -1.5}
      }
      \picscale{1 -1}{
           \picfillgraycol{0.3 0.8 0.7}
           \picfilledcircle{1.5 0}{0.8}{$P$}
      }
      \picfillgraycol{0.7 0.2 0.6}
      \picfilledcircle{-1.5 0}{0.8}{$Q$}
   }
}
\,,
\]
which are both mutations. (Altering the way of building connected sums
out of the prime factor diagrams is also considered a mutation.)

Since mutation does not alter the determinant, we have by
induction $\dt(D')\ge (n-1)(n-4)$.

Now let $D'_b$ and $D'_w$ we the diagrams obtained from $D$
when at both $p$ and $p'$ the black resp. white splicings
are applied, and $L_b$ and $L_w$ the alternating links they
represent.

When a nugatory white splicing at a crossing $r$ in $D$ renders a
crossing $q$ nugatory, then we have a white curve $\gm$ as
in \eqref{gmeq} though $p$ and $q$.
Similarly if two white splicings at crossings $r$ and $s$
render $q$ nugatory, then we have a white curve $\gm$ in $D$
through $q$, $r$ and $s$.

Assume now, a crossing $q\not\in\{p,p'\}$ is nugatory in both
$D'_b$ and $D'_w$. Then we have a white curve $\gm_{q,w}$ and
a black curve $\gm_{q,b}$ in $D$ through $q$ and at least one of
$p$ and $p'$ (and no other crossing different from $p$ and $p'$).
By the Jordan curve theorem, $|\gm_{q,w}\cap\gm_{q,b}|=2$.

If $\gm_{q,b}\cap D=\gm_{q,w}\cap D$ and $|\gm_{q,w}\cap D|=2$,
then $D$ is a Hopf link diagram, which we excluded.

Thus $\{|\gm_{q,b}\cap D|,|\gm_{q,w}\cap D|\}=\{2,3\}$. We claim
that for each of the two choices $|\gm_{q,b}\cap D|=3$ and
$|\gm_{q,w}\cap D|=3$ there is at most one $q\not\in\{p,p'\}$
satisfying these conditions. This is so, because $\gm_{q,b}\cap
D=\{q,p,p'\}$ and $\gm_{q',b}\cap D=\{q',p,p'\}$ imply that $q$ and
$q'$ have the same pair of opposite black regions. (These are the
regions adjacent to exactly one of $p$ and $p'$.) Then every
white curve through $q'$ must pass through $q$, contradicting
$\gm_{q',w}\cap D\subset \{q',p,p'\}$ if $q\ne q'$.
(For $|\gm_{q,w}\cap D|=3$ switch black and white.)

% 
% Since $D$ is not a Hopf link diagram, there is no similar
% curve $\gm'$ going through $q$ and only one further crossing and
% through black regions.
% 
% In a similar way, $\gm$ and $\gm'$ cannot be made to pass both
% though $q$ and both $p$ and $p'$ each. If $\gm$ passes though $q$,
% $p$ and $p'$, then 

Therefore, each $q\not\in\{p,p'\}$ is not nugatory at least in
one of $D'_b$ and $D'_w$, with at most 2 exceptions. Hence,
$c(L_b)+c(L_w)\ge 2n-4$, and by Crowell's result $\dt(D'_b)+
\dt(D'_w)\ge 2n-4$. 

Thus finally
\begin{myeqn}{\qed}
\dt(D)\,\ge\,\dt(D')+\dt(D'_b)+\dt(D'_w)\,\ge\,
(n-1)(n-4)+2n-4\,=\,n(n-3)\,.
\end{myeqn}

As a graph theoretical application, we obtain:

\begin{corr}
The number of spanning trees in a planar connected self-dual
graph with $2n$ edges is at least $n(n-3)$. \qed
\end{corr}

One can also reformulate theorem \reference{Thsq} graph-theoretically.
For this one remarks that by \cite{Menasco} an alternating diagram
of a composite alternating knot is composite, and checkerboard
graphs of such diagrams have a \em{cut vertex}. (This is a vertex,
which when removed together with all its incident edges, disconnects
the graph.) With the same argument as in the proof of theorem
\reference{th5.1} then one has:

\begin{theorem}\label{th5.1'}
Let $n$ be an odd perfect square. Then $n$ is the number of
spanning trees in a planar self-dual graph without cut vertex
if and only if $n\ne 1,\,9,\,49$. \qed
\end{theorem}

In the same way one can pose a conjecture which is
slightly stronger than conjecture \reference{Cjsq}:

\begin{conj}\label{Cjsq'}
Let $n$ be an odd natural number. Then $n$ is the number of
spanning trees in a planar self-dual graph without cut vertex
if and only if $n$ is the sum of two squares and $n\ne 1,\,9,\,49$.
\end{conj}

\section{The leading coefficients of the Alexander
and HOMFLY polynomial\label{S4}}

In this final section we explain another, apparently unrelated,
but also very striking occurrence of squares in connection with
invariants of achiral knots, namely in the leading coefficients of their
Alexander polynomial, and discuss a possible generalization of
this property to the HOMFLY, or skein, polynomial. Our results
have been obtained independently (but later) by C. Weber
and Q. H. C\^am V\^an \cite{VW}.

The skein polynomial $P$ is a Laurent polynomial in two variables $l$
and $m$ of oriented knots and links and can be defined
by being $1$ on the unknot and the (skein) relation
\begin{eqn}\label{1}
l^{-1}\,P\bigl(
\diag{5mm}{1}{1}{
\picmultivecline{0.18 1 -1.0 0}{1 0}{0 1}
\picmultivecline{0.18 1 -1.0 0}{0 0}{1 1}
}
\bigr)\,+\,
l \,P\bigl(
\diag{5mm}{1}{1}{
\picmultivecline{0.18 1 -1 0}{0 0}{1 1}
\picmultivecline{0.18 1 -1 0}{1 0}{0 1}
}
\bigr)\,=\,
-m\,P\bigl(
\diag{5mm}{1}{1}{
\piccirclevecarc{1.35 0.5}{0.7}{-230 -130}
\piccirclevecarc{-0.35 0.5}{0.7}{310 50}
}
\bigr)\,.
\end{eqn}
This convention uses the variables of \cite{LickMil}, but
differs from theirs by the interchange of $l$ and $l^{-1}$.
% We call the three diagram fragments in \eqref{1} from left to
% right a \em{positive} crossing, a \em{negative} crossing and a
% \em{smoothed} out crossing.

There is a classic substitution formula (see \cite{LickMil}),
expressing the Alexander polynomial $\Dl$, for the normalization
so that $\Dl(t)=\Dl(1/t)$ and $\Dl(1)=1$, as a special case of
the HOMFLY polynomial:
\begin{eqn}\label{PDl}
\Dl(t)\,=\,P(i,i(t^{1/2}-t^{-1/2}))\,.
\end{eqn}

We will denote by $\md_mP$ the \em{maximal degree} of $m$ in $P$,
and by $\mc_m P$ the \em{leading coefficient} of $m$ in $P$.
Similarly we will write $\md\Dl$ and $\mc\Dl$.

\subsection{Problems and partial solutions}

If $K$ is an alternating knot, then the HOMFLY polynomial $P_K\in
\bZ[m^2,l^{\pm 2}]$ is known to be of the form 
\[
a_{2g}(l)m^{2g}+(\mbox{lower $m$-degree terms})\,,
\]
with $a_{2g}\in\bZ[l^2,l^{-2}]$ being a non-zero
Laurent polynomial in $l^2$ and $g=g(K)$ the \em{genus} of $K$,
the minimal genus of an embedded oriented surface $S\subset \bR^3$ with
$\partial S=K$. (See \cite{Cromwell}.)
If $K$ is achiral, then $a_{2g}(l^{-2})=a_{2g}(l^2)$, that is,
$a_{2g}$ (and, in fact, all the other coefficients of $m$ in $P_K$)
is self-conjugate.

The main problems we consider here can be formulated as follows.

\begin{question}\label{ques1}
Is $\mc\Dl_K$ for an achiral knot $K$ always a square up to sign,
and if $\Dl$ is normalized so that $\Dl(t)=\Dl(1/t)$ and $\Dl(1)=1$,
is the sign $\sgn (\mc \Dl_K)$ always given by $(-1)^{\md \Dl_K}$?
\end{question}

The questions on $\Dl$ can be generalized to $P$.

\begin{question}\label{ques3}
For which large knot classes is it true that achiral knots
have $\mc_m P$ of the form $f(l^2)f(l^{-2})$ for some $f\in\bZ[l]$?
\end{question}

This is true for several special cases.
It appears convenient to compile them into one single statement.

Recall, that a knot $K$ is \em{fibered}, if $S^3\sm K$ fibers over $S^1$
(with fiber being a minimal genus Seifert surface for $K$), and
\em{homogeneous}, if it has a diagram $D$ containing in each connected
component of the complement (in $\bR^2$) of the Seifert circles of $D$
(called \em{block} in \cite[\S 1]{Cromwell}) only crossings of the
same sign.

\begin{propo}\label{op}
Let $K$ be an ($+/-$)achiral knot. Then $\mc_m P_K$ is of the form
$f(l^2)f(l^{-2})$ for some $f\in\bZ[x]$, if
% (achiral alternating)
{\nopagebreak
\def\labelenumi{\arabic{enumi})}\mbox{}\\[-18pt]
\def\theenumi{\arabic{enumi})}
\begin{enumerate}
\item\label{casei} $K$ is a fibered homogeneous knot,
\item\label{caseii} $K$ is a homogeneous knot of crossing number at most 16, or
% \item\label{caseiii} $K$ is a rational (or $2$-bridge) knot, or
\item\label{caseiv} $K$ is an alternating knot. %(with an alternating
%diagram) not admitting a flype, or equivalently, by \cite{MenThis},
%with only one alternating diagram
\end{enumerate}
}
\end{propo}

{}From formula \eqref{PDl} it is straightforward that whenever the
leading $m$-coefficient of $P_K$ is of the above form, both the
modulus and sign of $\mc \Dl_K$ are as requested in
question \reference{ques1}.
% we have
% that $|\mc \Dl_K|$ is a square, and if $\Dl$ is normalized as
% explained, that $\sgn(\mc \Dl_K)=(-1)^{\md \Dl_K}$.
For these properties we have some more situations where they
can be established. We give again the so far complete
list of such cases, even if some of them are trivial.
% so again we give them as a list.

\begin{propo}\label{oo}
Let $K$ be an achiral knot. Then $|\mc \Dl_K|$ is a square, if
{\nopagebreak
\def\labelenumi{\arabic{enumi})}\mbox{}\\[-18pt]
\def\theenumi{\arabic{enumi})}
\begin{enumerate}
\item\label{Casei} $K$ is a fibered knot,
\item\label{Caseii} $K$ is a knot of crossing number at most 16,
% \item\label{Caseiii} $K$ is a rational (or $2$-bridge) knot,
\item\label{Caseiv} $K$ is an alternating knot,
\item\label{Casev} $K$ is strongly achiral, or
\item\label{Casevi} $K$ is negative achiral.
\end{enumerate}
}
Moreover, $\sgn(\mc \Dl_K)=(-1)^{\md \Dl_K}$, if
{\nopagebreak
\def\labelenumi{\arabic{enumi})}\mbox{}\\[-18pt]
\def\theenumi{\arabic{enumi})}
\begin{enumerate}
\setcounter{enumi}{5}
\item\label{CCaseiii} $K$ is a fibered homogeneous knot,
\item\label{CCaseii} $K$ is a knot of crossing number at most 16,
\item\label{CCasei} $K$ is an alternating knot,
\item\label{CCaseiv} $K$ is strongly achiral, or 
\item\label{CCasev} $K$ is negative achiral. 
\end{enumerate}
}
\end{propo}

% Let $K$ be an ($+/-$)achiral knot. Then $|\mc \Dl_K|$ is a square, if
% 
% In particular, for any of these knots, and additionally for any
% (other) achiral knot $K$ of crossing number at most 16,
% $\mc \Dl_K\cdot (-1)^{\md \Dl_K}$ is a square (of an integer),
% where $\Dl$ is normalized so that $\Dl(t)=\Dl(1/t)$ and $\Dl(1)=1$.

\begin{rem}
We omitted to explicitly mention rational knots in proposition
\reference{oo}, as we already remarked that they are
alternating, and the achiral ones are strongly ($-$)achiral.
\end{rem}

In the following subsections we will collect the arguments
establishing the conjectured properties in the indicated
special cases. Some of them are well-known, or a matter of electronic
verification, and thus do not deserve separate proof. These parts
are briefly discussed first. Our main result are the statements in
the alternating case, which are proved subsequently.

\subsection{Some known and experimental results}

Question \reference{ques1} was the (chronologically) first question
I came across, addressing special properties of the leading
coefficients of the Alexander and skein polynomial of achiral knots.

This question came up when considering the formula
\[
\mc\Dl_{K}\,=\,\pm 2^{-2g}\,\prod_{i=1}^{2g}\,a_i\,
\]
for a rational knot $K=(a_1 \dots a_{2g})$ with all $a_i\ne 0$ even.
If $K$ is achiral, the sequence $(a_1,\dots,a_{2g})$ is palindromic,
and so we have, up to the sign, the requested property for rational
knots. A further larger class of achiral knots satisfying
the conjectured condition are the strongly achiral knots (see the
Theorem of \cite{HarKaw}). Subsequently, I verified all (prime) knots
in Thistlethwaite's tables \cite{KnotScape} up to 16 crossings
(note, that the property for a composite knot will follow from that
of its factors), and found no counterexample.

Although there seems much evidence for a positive answer to question
\reference{ques3}, its
diagrammatic, and not topological, origin (see \S\reference{S5})
suggests that it may not be true in general, but at least on some 
(diagrammatically defined) nice knot classes, for example alternating
knots.

We now collect the arguments that prove the easier cases of
propositions \reference{op} and \reference{oo}.

\proof[of proposition \reference{op} except part \reference{caseiv} and
proposition \reference{oo} except part \reference{Caseiv}]

% A further-going question is
% 
% \begin{question}\label{ques2}
% If for an achiral knot the answer to question \reference{ques1} is `yes',
% and  $\Dl$ is normalized so that $\Dl(t)=\Dl(1/t)$ and $\Dl(1)=1$,
% is $\sgn (\mc \Dl_K)=(-1)^{\md \Dl_K}$?
% \end{question}
% 
{\em alternating knots}. As said, we will deal with the squareness
and HOMFLY polynomial later. As for the further-going question on the
sign, the positive answer follows from the alternation of the
coefficients of $\Dl$ proved by Crowell in \cite{Crowell}, and
the property $\Dl(-1)>0$ following from Murasugi's trick,
as explained in \S\reference{S2}.

\begin{figure}[htb]
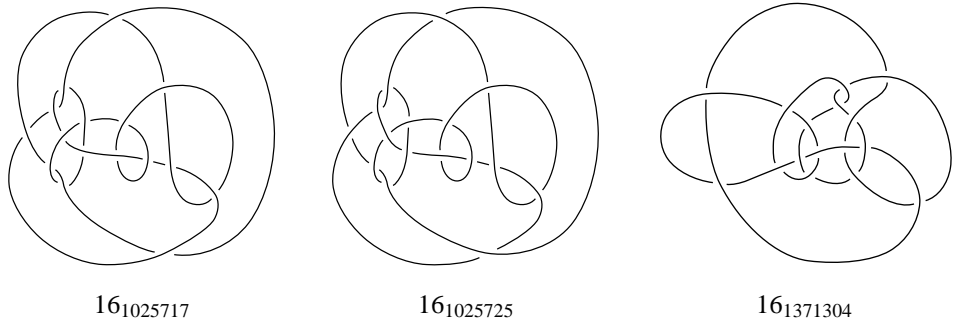

\[
\begin{array}{*2{c@{\qquad}}c}
\epsfsv{3.6cm}{t1-16_1025717} & 
\epsfsv{3.6cm}{t1-16_1025725} & 
\epsfsv{4cm}{t1-16_1371304} \\
\ry{6mm} 16_{1025717} & 16_{1025725} & 16_{1371304}
\end{array}
\]
% \tm
\caption{Three $+$achiral knots with $\mc_m P$ not of
the form $f(l^2)f(l^{-2})$.\label{fig1}}
\end{figure}

{\em knots with at most $16$ crossings}.
The answer is also `yes' for $\le 16$ crossing knots. This follows
from some experimental results related to question \reference{ques3}.
It is clear that an answer `yes' to question \reference{ques3}
implies the same answer to question \reference{ques1}.
% and \reference{ques2}.
This time a computer experiment found that the answer is not positive
in general, but the examples showing exceptional behaviour are not quite
simple, and required to use the full extent of the tables presently
available. Among $\le 16$ crossing knots, only three fail to have
this property: $16_{1025717}$, $16_{1025725}$ and $16_{1371304}$.
They are all $+$achiral and have $P=m^8(l^{-2}+3+l^2)+O(m^6)$.
See figure \reference{fig1}. Since the Alexander polynomial of
the three knots has degree 4 and leading coefficient $1$, they still
conform to the properties requested in question \reference{ques1}.
Also, at least the first two knots in figure \reference{fig1}
were found to be fibered by the method of \cite{Gabai}, showing
that the homogeneity assumption is essential in part \reference{casei}
of proposition \reference{op}.

{\em fibered knots}.
On the other hand, a class of knots where the squareness property
of $\mc\Dl$ is trivial, are the fibered knots, since then
$\mc\Dl=\pm 1$. For fibered homogeneous (in particular, alternating)
knots, the other properties also follow easily from known results,
because by \cite[corollaries 4.3 and 5.3]{Cromwell} and
\cite{MurPrz} we have for such knots that $\mc_m P=l^k$ for some
$k\in2\bZ$, and then achirality shows $k=0$. 

%Moreover, also for rational knots the answer is again `yes'
% (see below).

{\def\md{\min\deg}\def\Md{\max\deg}{\em strongly achiral knots}.
For strongly achiral knots the claims of proposition \reference{oo}
follow directly from the results of \cite{HarKaw} stated in theorem
\reference{thHK}. The only non-obvious property may be the sign of
$\mc\Dl$ for a strongly negative amphicheiral knot. To see this, first
normalize the polynomial $F$ found by the theorem by some $+t^n$
so that $\md F=-\Md F$. Then we have
\begin{eqn}\label{xxx}
F(t)=\pm F(-t^{-1})\,.
\end{eqn}
If we normalize $\Dl$ so that $\md\Dl=-\Md\Dl$ and $\Dl(1)=1$,
then the minimal and maximal degrees show that we must have $n=0$ in
\[
\Dl(t^2)\,=\,\pm t^n \,F(t)F(t^{-1})\,,
\]
and the value at $t=1$ shows that we must have the positive sign.
Denote by $[X]_i$ the coefficient of $t^i$ in $X\in\cL$.
Then, since $\Dl(1)=1$ and $\Dl(t)=\Dl(t^{-1})$, the absolute term 
$[\Dl(t)]_0$ of $\Dl(t)$ is odd. Thus the same is true for
$\Dl(t^2)=F(t)F(t^{-1})$. But
\[
\bigl[F(t)F(t^{-1})\bigr]_0\,=\,\sum_{i=\md F}^{\Md F}[F(t)]_i^2\,,
\]
and so from \eqref{xxx} we conclude that 
$F$ must have non-zero absolute term. This determines the sign
in \eqref{xxx} to be positive, and then $\mc F=\pm \mcf F$
dependingly on the parity of $\Md F=\Md\Dl$.

{\em negative amphicheiral knots}. Hartley \cite{Hartley} has
extended the result of \cite{HarKaw} for strongly negative amphicheiral
knots to arbitrary negative amphicheiral knots. Thus the claim
follows from the previous argument. \qed

\begin{rem}
In fact, in \cite{Kaw}, Kawauchi conjectures that the property of
the Alexander polynomial of a strongly negative amphicheiral
knot he proves with Hartley in \cite{HarKaw}, and later Hartley
\cite{Hartley} generalizes to an arbitrary negative amphicheiral knot,
extends to the Alexander polynomial of an arbitrary amphicheiral knot.
This conjecture clearly implies a positive answer to question
\reference{ques1}. Kawauchi's conjecture is true in particular for
2-bridge knots, since in \cite{HarKaw} he shows that all amphicheiral
2-bridge knots are strongly negative amphicheiral. I verified the
conjecture for all prime amphicheiral knots of $\le 16$ crossings.
Note that Hartley also obtains a condition for positive amphicheiral
knots, but it is too weak to address any of our questions.
\end{rem}
}

\begin{rem}
Fibered homogeneous knots
contain the homogeneous braid knots of \cite{Stallings}, but also
many more. For example, there are 15 fibered homogeneous prime 10
crossing knots, among them 12 alternating and 2 positive ones, which
can be shown by the work of \cite{Cromwell} and an easy computer check
not to have homogeneous braid representations: $10_{60}$, $10_{69}$,
$10_{73}$, $10_{75}$, $10_{78}$, $10_{81}$, $10_{89}$, $10_{96}$,
$10_{105}$, $10_{107}$, $10_{110}$, $10_{115}$, $10_{154}$, $10_{156}$
and $10_{161}$.
\end{rem}

\subsection{Two general statements\label{TGS}}

The remaining cases of propositions \reference{oo} and
\reference{op} are included in two more general theorems.
The first one generalizes the result of \cite{MurPrz2}, where
it was shown that for an alternating amphicheiral knot 
the leading coefficient of the Alexander polynomial is not a prime.

\begin{theorem}\label{th6.1}
Let a knot $K$ have a homogeneous diagram $D$ which can be turned into
its mirror image (possibly with opposite orientation) by a sequence of
mutations and $S^2$-moves (changes of the unbounded region).
Then $|\mc \Dl_K|$ is a square.
\end{theorem}

\proof[of theorem \reference{th6.1}]
Let $\tl D$ denote the mutation equivalence class of a knot
diagram $D$ (that is, the set of all diagrams that can be obtained
from $D$ by a sequence of mutations). Assume $D$ to be non-split
(the split case easily reduces to the non-split one). We consider
orientation reversal as a special type of mutation, so a diagram
and its inverse belong to the same mutation equivalence class.

As in \cite[\S 1]{Cromwell}, the Seifert picture of $D$ defines a
decomposition of $D$ into the \em{$*$--product (or Murasugi sum)}
of special alternating
diagrams $D_1,\dots,D_n$, called \em{blocks}. These diagrams may not be
prime. Let $D_{i,1},\dots,D_{i,n_i}$ be the prime components of $D_i$.
Note that all $D_{i,j}$ are positive or negative (dependingly
on $D_i$). They will have no nugatory crossings if $D$ has neither.

Define
\[
\cI(D)\,:=\,\{\,\tl D_{i,j}\,\}_{i=1,\dots,n,\,j=1,\dots,n_i}\,.
\]
Here a set is to be understood with the order of its elements
ignored, but with their multiplicity counted (i.e., $\{1,1,2,3\}=
\{1,1,3,2\}\ne \{1,2,3\}$).

Now apply a mutation on $D$. The Seifert picture separates
the Conway circle into 3 parts $A$, $B$ and $C$.
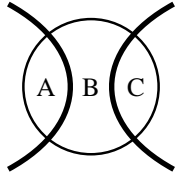
\[
\diag{1cm}{2}{2}{
  \pictranslate{1 1}{
    \piccircle{0 0}{0.9}{}
    \piclinewidth{14}
    \picmultigraphics[S]{2}{-1 1}{
       \piccurve{1.1 d}{0. 0.5}{0. -0.5}{1.1 -1.1}
    }
  }
  \picputtext{1.6 1}{C}
  \picputtext{0.4 1}{A}
  \picputtext{1 1}{B}
}
\]
Because the Conway circle intersects the Seifert picture only
in 4 points, all parts $A$, $B$ and $C$ represent connected
components of the blocks in $D$ they belong to (or possibly
connected sums of several such connected components).

Mutation then has the effect of applying mutation on $B$ and
interchanging and/or reversing $A$ and $C$. Therefore, $\cI(D)=
\cI(D')$ for any iterated mutant diagram $D'$ of $D$.

If $D$ has the property assumed in the theorem, then $\cI(D)=
\cI(!D)$, or
\[
\{\,\tl D_{i,j}\,\}_{i=1,\dots,n,\,j=1,\dots,n_i}\,=\,
\{\,\wt{!D_{i,j}}\,\}_{i=1,\dots,n,\,j=1,\dots,n_i}\,.
\]
Let $\phi\,:\,\{\,\tl D_{i,j}\,\}\to \{\,\tl D_{i,j}\,\}$
be the bijection induced by $\tl D_{i,j}\mapsto \wt{!D_{i,j}}$.

Since mutation preserves the writhe, $\phi$ has no
fixpoints (unless some $D_{i,j}$ has no crossings, in which
case $D$ is split). Thus $\phi$ descends to a bijection
\[
\phi\,:\,\{\,\tl D_{i,j}\,:\,\mbox{$D_{i,j}$ positive}\,\}\to 
\{\,\tl D_{i,j}\,:\,\mbox{$D_{i,j}$ negative}\,\}\,.
\]
Then by \cite{Murasugi2}, $\mc\Dl$ is multiplicative under
$*$-product, and hence
\begin{eqnarray*}
\mc \Dl_D & = & \prod_{i,j}\,\mc \Dl_{D_{i,j}}\\
          & = & \prod_{i,j\,:\,\scbox{$D_{i,j}$ positive}}
\mc \Dl_{D_{i,j}}\,\cdot\,\mc \Dl_{!D_{i,j}}\\
          & = & \prod_{i,j\,:\,\scbox{$D_{i,j}$ positive}}
\mc \Dl_{D_{i,j}}\,\cdot\,\pm \mc \Dl_{D_{i,j}}\\
          & = & \pm\left(\,\prod_{i,j\,:\,\scbox{$D_{i,j}$ positive}}\,
\mc \Dl_{D_{i,j}}\right)^2\,,
\end{eqnarray*}
as desired. \qed

\begin{theorem}\label{thP}
Under the same assumption as theorem \reference{th6.1} we have
$\mc_m P(D)=f(l^2)f(l^{-2})$ for some $f\in\bZ[x]$.
\end{theorem}

\proof Using \cite{MurPrz2} instead of \cite{Murasugi2},
we obtain $\mc_m P(D)=f(l)f(l^{-1})$. Since $\mc_m P(D)$
has only even powers of $l$, the result follows. \qed

\begin{corr}\label{crx}
For any alternating achiral knot $K$, %we have that
$|\mc \Dl_K|$ is a square and $\mc_m P_K=f(l^2)f(l^{-2})$.
\end{corr}

\proof Use that an alternating diagram is homogeneous, \cite{MenThis},
and that a flype can be realized as a sequence of mutations. \qed

Note, that in the case of $\Dl$ we obtain an exact condition
when a number occurs as $|\mc \Dl_K|$ for an alternating achiral knot
$K$ (since the other implication is trivial).

\begin{rem}
We need the homogeneity of $D$ only to assure that all
$D_{i,j}$ are positive or negative. (This weaker property is
invariant under mutations, whereas homogeneity is not.)
Thus the theorems could be formulated even slightly more
generally, but then also more technically.
\end{rem}

Although the corollary is the most interesting special case of
the theorems, they give indeed more general statements.

\begin{exam}
The non-alternating achiral knots $14_{45317}$ and $14_{45601}$ have
unique minimal diagrams (which therefore must be transformable into
their mirror
images by $S^2$-moves only), which are both homogeneous (of genus
$5$ and $4$, respectively).
\end{exam}

\autoepsftrue
\begin{figure}[htb]
\[
\begin{array}{*2{c@{\qquad}}c}
\epsfsv{4cm}{t1-14_45317} & 
\epsfsv{4cm}{t1-14_45601} & 
\epsfsv{4cm}{t1-14_41330} \\
\ry{6mm} 14_{45317} & 14_{45601} & 14_{41330}
\end{array}
\]
% \tm
\caption{\label{fig2}}
\end{figure}

\begin{rem}
Here we consider (and mean unique) minimal diagrams only up to
$S^2$-moves, and not as in \cite{KnotScape} up to $S^2$-moves and
mirroring. There are achiral knots with a unique minimal diagram
up to $S^2$-moves and mirroring, but corresponding to two different
(mirrored) diagrams up to $S^2$-moves only, which are not
interconvertible by flypes. One such example is $14_{41330}$.
\end{rem}

\begin{rem}\label{remfl}
If a knot has a diagram $D$, which can be transformed into its
(possibly reverted) obverse by moves in $S^2$ and flypes, then
it also has a diagram $D'$, which can be transformed into
its obverse by moves in $S^2$ only. This follows from the fact
that mirroring and moves in $S^2$ take the flyping circuits
of $D$ into each other (for the definition of flyping circuits
see \cite[\S 3]{SunThis}), and flyping in a flyping circuit
is independent from the other ones. $D'$ can be obtained by
appropriate flypes from $D$. (For an analogous statement
about the checkerboard graphs, see \cite{DH}.)
\end{rem}

\subsection{Some diagrammatic questions\label{S5}}

The two theorems of \S\reference{TGS} suggest the
diagrammatic arguments motivating questions \reference{ques1} and
\reference{ques3}. We conclude with a more detailed problem
concerning possible generalizations. In order to make the result
of \cite{MurPrz} work, we need to consider $P$-maximal diagrams.

\begin{defi}
Call a link diagram $D$ with $c(D)$ crossings and $s(D)$ Seifert circles
\em{$P$-maximal}, if $\md_mP(D)=c(D)-s(D)+1$.
\end{defi}

In \cite{Morton}, Morton showed that the one inequality $\md_mP(D)\le
c(D)-s(D)+1$ holds for any arbitrary link diagram, and used this
to show that there are knots $K$, which do not possess a diagram $D$
with $g(D)=g(K)$ (a fact that also implicitly follows from
\cite{Whitten}). Here $g$ denotes the \em{genus} of a knot or knot
diagram, for latter being defined by %where for a \em{knot} diagram $D$,
$g(D)=\myfrac{1}{2}(c(D)-s(D)+1)$, which is the genus of the
canonical Seifert surface associated to $D$ (see \cite[\S 4.3]{Adams}
or \cite{Rolfsen}).

In \cite{Cromwell} it was shown that homogeneous diagrams are
$P$-maximal. Many knots have $P$-maximal diagrams~-- beside the
homogeneous knots, for example all (other) knots in Rolfsen's tables
\cite[appendix]{Rolfsen} and also all the 11 and 12 crossing knots
tabulated in \cite{KnotScape}.
However, some knots do not~-- in \cite [fig.~9]{gen2} we gave four
examples of 15 crossings.

In \cite{MurPrz} it was shown that $\mc_mP$ is multiplicative under
$*$-product of $P$-maximal diagrams. This is the link
between the above polynomial conjectures and the diagrammatic problems,
which we summarize in the question below.

\begin{question}
Does any ($+/-$)achiral knot (or an achiral knot in which large
knot class) have a diagram that can be
{\nopagebreak
\def\labelenumi{\Alph{enumi})}\mbox{}\\[-10pt]
\def\theenumi{\Alph{enumi})}
\def\labelenumii{\arabic{enumii})}
\def\theenumii{\arabic{enumii})}
\begin{enumerate}
\item\label{case1} transformed into its (possibly reverted) obverse by
moves in $S^2$ (changes of the unbounded region), or

% This question is mainly motivated by \cite{MenThis}, where a positive
% answer was given for alternating knots. It is interesting whether
% the arguments after \reference{case1} can be extended to show that
% a positive answer for homogeneous (or at least alternating, where
% it was already given) diagrams implies a positive answer to question
% \reference{ques3} (note, that flypes in general do not preserve
% homogenuity). As seen, there is strong evidence for question
% \reference{ques3} to have a positive answer for alternating knots.
% \end{enumerate}
\item\label{caseB} represented as the iterated connected and
Murasugi sum of (some even
number of) $P$-maximal link diagrams $D_i$ with $\{D_i\}$ being
mutually obverse (up to orientation and $S^2$-moves) in pairs?
\end{enumerate}
}
\end{question}

As a motivation, we mention briefly some % problem we discuss
related partial cases and implications. % but also some cautionary
% evidence.

Remarks on part \reference{case1}:
\begin{itemize}
\item By considering the blocks of such a diagram,
we see that it is the $*$-product of special diagrams
$D_i$, such that each $D_i$ is transformable by $S^2$-moves into the
obverse of itself, or of some other $D_j$, $j\ne i$.
\item If the answer is `yes' for some \em{homogeneous} diagram (which
in particular happens by \cite{MenThis} for alternating diagrams not
admitting flypes), no blocks transform onto their own obverses, and
\cite{MurPrz} shows a positive answer to question \reference{ques3}.

The $(l-)$coefficients of $\mc_mP$ of the three knots in figure
\reference{fig1} do not alternate in sign, and so these knots
cannot be homogeneous (beware of the different convention for $P$ in
\cite{Cromwell}!). Still, the first two knots have a $P$-maximal
diagram of the type requested in part \reference{case1}. It can be
obtained from the one given in the figure by a flype. Thus theorem
\reference{thP} does not hold under the weaker assumption
the diagram to be $P$-maximal instead of homogeneous.
% but still provide some evidence against a general positive
% answer to this question.
\item The answer is `yes' for rational knots.
The (palindromic) expression $(a_1\dots a_na_n\dots a_1)$ with all
$2g$ even numbers
$a_i\ne 0$ gives a rational diagram (of the form $D*!D$, where $D$ is
a connected sum of diagrams of reversely oriented $(2,a_{2i})$-torus
links), having the desired property.

\item By remark \reference{remfl}, one can equivalently also allow
moves in $S^2$ and flypes, so the answer is in particular positive
for alternating knots by \cite{MenThis}.
\end{itemize}

Remarks on part \reference{caseB}:
\begin{itemize}
\item By \cite{MurPrz}, a positive answer implies a positive answer
also for question \reference{ques3} (thus in particular the answer
is negative for the knots on figure \reference{fig1}). This
is the motivation for proposing question \reference{ques3} after
arriving to question \reference{ques1}. 
\item In turn, a positive answer to \reference{caseB} is implied by
a positive answer to \reference{case1} for homogeneous diagrams.
This is the motivation for proposing part \reference{caseB}.
% \item The diagrams answering positively \reference{case1} for
% rational knots, answer positively \reference{caseB}, and hence also
% question \reference{ques3} in this case.
\end{itemize}

\noindent{\bf Acknowledgements.} I would like to thank to Kunio
Murasugi, Morwen Thistlethwaite, Kenneth Williams, and especially
to Don Zagier for their helpful remarks.
% for some helpful discussions.

\end{document}